\documentclass[onefignum,onetabnum]{siamart251216}

\usepackage{lipsum}
\usepackage{amsfonts}
\usepackage{graphicx}
\usepackage{epstopdf}
\usepackage{algorithmic}
\usepackage{makecell} 
\ifpdf
  \DeclareGraphicsExtensions{.eps,.pdf,.png,.jpg}
\else
  \DeclareGraphicsExtensions{.eps}
\fi

\newcommand{\omap}{\mathcal{G}^\dagger}
\newcommand{\cA}{\mathcal{A}}
\newcommand{\cU}{\mathcal{U}}

\newcommand{\R}{\mathbb{R}}

\newcommand{\E}{\mathbb{E}}

\newcommand{\cP}{\mathcal{P}}
\newcommand{\cR}{\mathcal{R}}
\newcommand{\cQ}{\mathcal{Q}}
\newcommand{\cL}{\mathcal{L}}
\newcommand{\cJ}{\mathcal{J}}

\newcommand{\cF}{\mathcal{F}}

\newcommand{\cG}{\mathcal{G}}

\newcommand{\bigO}{\mathcal{O}}
\newcommand{\dd}{\mathrm{d}}

\newsiamremark{remark}{Remark}
\newsiamremark{hypothesis}{Hypothesis}
\crefname{hypothesis}{Hypothesis}{Hypotheses}
\newsiamthm{claim}{Claim}

\headers{Cost-Accuracy Trade-offs}{D. Z. Huang, A. M. Stuart}

\title{Cost-Accuracy Trade-offs: Neural Operator vs Classical Numerical Solver\thanks{Submitted to the editors DATE.
\funding{
DZH acknowledges funding support from National Key R\&D Program of China 2025YFA1018700, and National Natural Science Foundation of China No.62595771.
AMS acknowledges funding support from a USA Department of Defense Vannevar Bush Faculty Fellowship (award N00014-22-1-2790). This work was initiated when DZH was a postdoctoral scientist at Caltech, also supported by award N00014-22-1-2790.
}}}

\author{Daniel Zhengyu Huang\thanks{Beijing International Center for Mathematical Research,  Center for Machine Learning Research, Peking University, Beijing, China (\email{huangdz@bicmr.pku.edu.cn}).}
\and Andrew M. Stuart\thanks{California Institute of Technology, Pasadena, CA 
  (\email{astuart@caltech.edu}).}}

\usepackage{amsopn}

\makeatletter
\newcommand*{\addFileDependency}[1]{
  \typeout{(#1)}
  \@addtofilelist{#1}
  \IfFileExists{#1}{}{\typeout{No file #1.}}
}
\makeatother

\ifpdf
\hypersetup{
  pdftitle={Cost-Accuracy Trade-offs: Neural Operator vs Classical Numerical Solver},
  pdfauthor={D. Z. Huang and A. M. Stuart}
}
\fi

\definecolor{darkred}{rgb}{.7,0,0}

\definecolor{darkblue}{rgb}{0,0,.7}

\begin{document}

\maketitle

\begin{abstract}
Neural operators are data-driven models that learn mappings from inputs that parameterize partial differential equations (PDEs), such as spatially varying coefficients, initial conditions, forcing terms, boundary conditions, or geometries, to solution fields or quantities of interest. Once trained, they can serve as surrogates for classical numerical solvers in many-query settings that require repeated evaluations for varying inputs. We address the question of when, and then why, neural
operator surrogates outperform classical numerical solvers, in terms of cost for a given accuracy.
We focus on the post-training, many-query limit, in which data-acquisition and training costs are treated as fixed and fully amortized.
Even in this deliberately favorable regime for neural operators, there are regimes in which
classical solvers outperform the surrogate models.
We compare the cost-accuracy performance of neural operator surrogates and classical numerical solvers through a reproducible benchmark study comparing neural operators with problem-matched classical solvers on representative problems in computational science and engineering, focusing on prediction error, per-query floating-point cost, and wall-clock runtime.
The results reveal no uniformly superior method.
Neural operators are most competitive at low-to-moderate accuracy requirements.
Their floating-point cost advantage depends strongly on the problem structure, arising when they avoid temporal or nonlinear iterations or predict a reduced quantity of interest rather than a full solution field. Additional wall-clock speedups result from dense tensor operations that are well suited to modern hardware. As the target accuracy is tightened, achieving the required accuracy with neural operators becomes increasingly challenging, and classical solvers outperform surrogates in this
regime; thus classical solvers will remain important for verification and high-accuracy computation,
even as surrogates emerge as an attractive alternative for parts of the engineering pipeline.

\end{abstract}

\begin{keywords}
Surrogate modeling, Neural operator, Parametric partial differential equations, Numerical methods
\end{keywords}

\begin{AMS}
65Y20, 68T07, 65N55, 65M70
\end{AMS}

\section{Introduction}
Many computational tasks in science and engineering require repeated evaluations of expensive computational models, often involving parametric partial differential equations. Representative examples include design optimization~\cite{li2022machine,zhou2024ai,luo2025efficient,gomes2026combining}, where many candidate designs must be screened efficiently; uncertainty quantification~\cite{zhu2018bayesian,cao2023residual,gao2024adaptive}, where probability estimates require a large number of model evaluations; and weather forecasting~\cite{kurth2023fourcastnet,bi2023accurate}, where predictions must be generated rapidly for evolving physical states, and ensembles of forecasts are often used. In such many-query settings, the dominant computational cost may arise not from a single high-fidelity simulation, but from repeatedly solving a related class of PDE problems for different physical states, parameters, or geometries. This motivates the construction of surrogate models to potentially replace classical numerical solvers.

Classical approaches to surrogate modeling include projection-based reduced-order
models~\cite{devore2014theoretical,carlberg2013gnat,peherstorfer2016data,binev2017data,cohen2020state,qian2020lift,grimberg2021mesh,gosea2021data}, Koopman-operator-based methods~\cite{schmid2010dynamic,rowley2009spectral,mezic2005spectral,mezic2013analysis}, and Gaussian-process emulators~\cite{kennedy2001bayesian,cleary2021calibrate,chen2021solving}.
More recently, neural-network-based surrogates~\cite{zhu2018bayesian,khoo2019switchnet,hesthaven2018non} have attracted
considerable attention because of their ability to approximate high-dimensional maps
and their compatibility with modern accelerator hardware.
\emph{Neural operators}~\cite{bhattacharya2021model,li2020fourier,kovachki2023neural,lu2021learning} extend this
perspective by learning mappings between input and output functions rather than between fixed finite dimensional vectors. Representative architectures include DeepONet~\cite{lu2021learning},  PCA-Net~\cite{hesthaven2018non,bhattacharya2021model}, graph-based neural operators~\cite{li2020neural,mousavi2025rigno},  kernel-based neural operators~\cite{li2020fourier,kovachki2023neural,li2022fourier,gin2020deepgreen,boulle2022learning,hao2026multiscale}, and attention-based neural operators~\cite{cao2021choose,wu2024transolver,calvello2025continuum}.
A central attraction of neural operators is that they formulate the learning problem as a mapping between function spaces and introduce discretization only for numerical computation, thereby providing a degree of robustness across discretizations and resolutions. 
However, this property alone does not establish a computational advantage over classical numerical solvers.

Training neural operators to a prescribed accuracy may require substantial investments in data generation, model capacity, and offline optimization. 
We deliberately set these offline costs aside and study the post-training, many-query limit, in which the dataset and trained model are fixed and query-independent costs are fully amortized. Our comparison between neural operator surrogates
and classical numerical solvers therefore concerns the cost of answering a new query at comparable accuracy.
Even in this deliberately favorable regime, recent studies have shown that some reported inference-efficiency advantages of neural PDE surrogates rely on comparisons with weak baselines or with non-comparable accuracy requirements~\cite{mcgreivy2024weak}. 
Cost-accuracy analysis is therefore essential
for understanding the practical value of these methods and guiding their use.
For classical numerical solvers, cost-accuracy behavior is relatively well understood
through a conjunction of sharp error estimates and computational cost estimates driven, primarily, by the cost of standard linear algebra and optimization routines. For neural operators, by contrast, a comparable understanding remains incomplete,
primarily because sharp error estimates are not available.
Existing theory primarily addresses approximation properties~\cite{chen1995universal,lu2021learning,kovachki2023neural,kovachki2021universal,lanthaler2022error,weihs2025deep} and statistical convergence rates~\cite{de2023convergence,liu2024deep,reinhardt2026statistical}, but does not directly characterize inference-time cost-accuracy trade-offs.
The present work provides a systematic numerical study of these cost-accuracy trade-offs, and addresses the following question: 

\begin{quote}
\textit{In the post-training, many-query limit, under what accuracy requirements and problem structures does neural operator inference achieve a favorable per-query cost-accuracy trade-off relative to a well-designed classical numerical solver?}
\end{quote} 

\noindent 
We measure per-query computational cost using floating-point work and CPU and GPU wall-clock runtime. Floating-point work characterizes algorithmic complexity, whereas wall-clock runtime also reflects implementation, parallelism, memory access, and hardware utilization.

\subsection{Contributions}
The main contributions and empirical findings are as follows:
\begin{itemize}
\item We establish three benchmarks for evaluating the relative merits of neural operators and classical numerical methods.
    \item At low-to-moderate accuracy requirements, neural operators can substantially outperform classical numerical solvers, in terms of wall-clock runtime at comparable error levels, in all three benchmarks.
    \item 
    At higher accuracy requirements, the
    tested neural operators are not cost-competitive with the classical numerical solvers in all three benchmarks.
\end{itemize}

More specifically, we formulate a post-training, many-query framework and apply it to three numerical experiments: 2D Darcy flow, 2D time-dependent incompressible Navier-Stokes flow, and 3D vehicle aerodynamics. In the first two experiments the
output is the entire solution field; in the vehicle case, the output is the surface pressure coefficient. The corresponding classical baselines are geometric multigrid, a Fourier pseudospectral solver, and an \texttt{OpenFOAM}-based finite volume solver. Across all three benchmarks, we use a modified Fourier neural operator architecture, enabling us to better isolate the effects of problem structure from those of architectural choice. The code, data generation procedures, solver baselines, neural operator architectures, and training protocols are available at
\begin{center}
\url{https://github.com/Zhengyu-Huang/Cost-accuracy-trade-off}, 
\end{center}
complementing earlier reproducible benchmark and software efforts~\cite{aristoff2023benchmark,lu2021deepxde,zou2024neuraluq,chen2025due}.

The observed online advantages arise through different computational mechanisms. For Darcy flow, the principal runtime benefit arises from hardware-efficient dense tensor operations rather than reduced floating-point work. For Navier-Stokes flow, inference reduces sequential depth by replacing many time steps with learned updates. For vehicle aerodynamics, direct prediction of the surface quantity of interest avoids a full volumetric solution and nonlinear iteration. These results show that the many-query limit alone does not guarantee an advantage; the outcome depends on which components of the classical computation are bypassed by the surrogate neural
operator.

The experiments also expose limitations that must be addressed before reliable deployment. For the tested neural operators, additional training data or model capacity yields diminishing improvements and may still fail to meet stringent accuracy requirements. In time-dependent problems, errors can accumulate rapidly during long-horizon prediction. More generally, performance may vary substantially across inputs, raising even greater concerns for out-of-distribution inputs; ignoring offline
training costs, as we have done here, is only reasonable if the trained model
generalizes well. Reliable use in many-query engineering workflows will 
therefore require uncertainty quantification~\cite{zou2024neuraluq,walz2024easy,winovich2026active}, methods for detecting inaccurate predictions, and hybrid strategies~\cite{zhang2024blending} that use data assimilation~\cite{revach2022kalmannet} and classical solvers to verify, correct, or replace neural predictions when higher accuracy and reliability are required.

\subsection{Organization}
The remainder of the paper is organized as follows. In Section~\ref{sec:PF}, we define
the many-query problem and the cost-accuracy comparison scope. In Section~\ref{sec:exp}, we outline the experimental setup, including the benchmarks, classical solvers, and neural operators. Section~\ref{sec:num_study} presents the numerical results, the corresponding cost-accuracy trade-offs, and underlying acceleration mechanisms and limitations of neural operators. We conclude in Section~\ref{sec:Conclusion}.

\section{Problem formulation}
\label{sec:PF}
In this section, we formulate the many-query PDE task and the comparison framework used throughout the paper. We first specify the parametric PDE and the requested output, and then introduce cost models for both classical solvers and neural operators. Finally, we consider the many-query limit and study the cost-accuracy 
trade-offs between neural operator surrogates and classical numerical solvers. We measure accuracy using a norm on the output space and computational cost in terms of floating-point work or wall-clock time.

\subsection{Many-query problem}
Consider a  parametric PDE together with an associated output map:
\begin{equation}
\label{eq:abstract-pde}
    \cR_D(u;a)=0,
    \qquad
    u_{\rm q}= \Phi_D(u).
\end{equation}
Here, $D\subseteq\mathbb R^d$ denotes the geometric configuration and determines a physical PDE domain $\Omega_D$, on which the PDE solution $u$ is defined. The input \(a\) collects the query-dependent problem data, which may, for example, be specified throughout $\Omega_D$ or on part of its boundary. The requested output \(u_{\rm q}\) may likewise be defined
throughout \(\Omega_D\), on part of its boundary, or as a finite dimensional vector. Examples of input $a$ include spatially varying coefficients, initial conditions, forcing terms, boundary conditions, or other problem parameters. When the geometry varies, $D$ is also part of the query. The output map \(\Phi_D\) extracts a prescribed quantity of interest and may simply be the identity map.

Solving \eqref{eq:abstract-pde} and evaluating $\Phi_D$ defines the input-output map
\begin{equation}
\label{eq:target-operator}
    \omap:(a,D)\longmapsto u_{\rm q}.
\end{equation}
In the benchmarks considered here, we assume that, after any problem-specific encoding, both the input \(a\) and output \(u_{\rm q}\) are represented as functions on the same domain \(D\). For the full-field
benchmarks, \(D=\Omega_D\). We define the input and output spaces by 
\begin{align*}
    \cA = \bigsqcup_{D\in\mathfrak D}  \bigl\{a: D \rightarrow \R^{d_a}\bigr\} \times \bigl\{D \bigr\}, \qquad \cU = \bigsqcup_{D\in\mathfrak D} \bigl\{u_{\rm q}: D \rightarrow \R^{d_u} \bigr\},
\end{align*}
where \(\mathfrak D\) is a family of admissible geometric configurations.
Let \(\mu\) be a probability measure on \(\cA\). The input-output pairs then have the joint law
\(\bigl(\mathrm{Id},\omap\bigr)^\#\mu\) on
\(\cA\times\cU\). The goal of surrogate modeling is to learn an efficient approximation to the target operator \(\omap\) from samples drawn from this law.

\subsection{Cost-accuracy model}
Suppose that a computational workflow must answer \(N_{\rm q}\) queries drawn from $\mu$ over its deployment lifetime, with an error tolerance \(\varepsilon\). The following cost model applies to both floating-point work and wall-clock runtime.

For each query, a classical numerical method discretizes and
solves the governing equations and then applies $\Phi_D$.  Its lifecycle cost is modeled as
\begin{equation}
\label{eq:lifecycle_solver}
\begin{aligned}
    C_{\rm solver}^{\rm life}(N_{\rm q}, \varepsilon)
    &= N_{\rm q}C_{\rm solver}(\varepsilon),
\end{aligned}
\end{equation}
where \(C_{\rm solver}\) is the mean, under $\mu$, cost per query required to attain the prescribed accuracy \(\varepsilon\).

A neural operator workflow additionally incurs offline costs for data acquisition and training. We decompose the lifecycle cost as
\begin{align}
C_{\rm NO}^{\rm life}(N_{\rm q},\varepsilon,N)
&=
C_{\rm data}(N, \varepsilon)
+
C_{\rm NO}^{\rm train}(N, \varepsilon)
+
N_{\rm q}C_{\rm NO}^{\rm infer}(\varepsilon),
\label{eq:lifecycle_NO}
\end{align}
where $C_{\rm data}$ is the cost of acquiring $N$ training pairs,
$C_{\rm NO}^{\rm train}$ is the optimization cost, and
$C_{\rm NO}^{\rm infer}$ is the mean cost of one inference. If all training outputs
are newly generated by the classical numerical solver at tolerance \(\varepsilon\), then
$C_{\rm data}(N,\varepsilon)\approx N C_{\rm solver}(\varepsilon)$. Existing or shared data lead to a different offline cost.

Multifidelity~\cite{peherstorfer2018survey} and
hybrid~\cite{zhang2024blending,hu2025hybrid} strategies are also used
to further improve the lifecycle cost-accuracy trade-off for neural operator
surrogates. However, we focus on a direct comparison between standalone classical solvers and standalone neural operators; coupled approaches are left for future work.

\subsection{Cost-accuracy comparison in the many-query limit}
It is useful to distinguish two deployment settings. In a data-scarce setting, a surrogate may be constructed to 
answer a finite number of queries, in which case data acquisition and training remain important
components of the lifecycle cost.
By contrast, a task-specific surrogate trained on a large dataset may be reused persistently for a much larger number of downstream queries, as in a digital wind tunnel used repeatedly for design screening, optimization, and uncertainty quantification. We focus on the idealized limit of the latter setting, $N_{\rm q}\rightarrow\infty$.

For fixed \(N\) and \(\varepsilon\), the per-query cost of the classical numerical solver is
\(C_{\rm solver}(\varepsilon)\). In contrast, the neural operator per-query cost is
\begin{align*}
\frac{C_{\rm NO}^{\rm life}(N_{\rm q},\varepsilon,N)}{N_{\rm q}}
&=
\frac{C_{\rm data}(N,\varepsilon)
+C_{\rm NO}^{\rm train}(N,\varepsilon)}
{N_{\rm q}}
+C_{\rm NO}^{\rm infer}(\varepsilon)\\
&\rightarrow
C_{\rm NO}^{\rm infer}(\varepsilon)
\quad\text{as }N_{\rm q}\to\infty.
\end{align*}
Thus, the per-query contribution of the fixed offline costs tends to zero. This limit deliberately favors the
neural-operator workflow by allowing data acquisition and training to be
fully amortized.

Therefore, in the many-query limit,  neural operator inference has a cost-accuracy advantage at
tolerance \(\varepsilon\),  when
\begin{align*}
\lim_{N_{\rm q} \rightarrow \infty}\frac{C_{\rm NO}^{\rm life}(N_{\rm q},\varepsilon, N) }{C_{\rm solver}^{\rm life}(N_{\rm q},\varepsilon)} < 1, 
\end{align*}
which is equivalent to
\begin{align}
\label{eq:NO_beat_solver}
C_{\rm NO}^{\rm infer}(\varepsilon) < C_{\rm solver}(\varepsilon).
\end{align}
The present work therefore asks whether \cref{eq:NO_beat_solver} holds, in which
accuracy regimes, and through which computational mechanisms.

\section{Experimental setup}
\label{sec:exp}
To test the many-query criterion in \cref{eq:NO_beat_solver}, we compare problem-matched classical solvers with the same kernel-based neural operator family on three PDE benchmarks. For each method, we vary the numerical or architectural configuration and report every resulting error-cost pair. Floating-point cost and wall-clock runtime cost, versus error, are reported separately.

\subsection{Three benchmarks and problem-matched baselines}

We consider three benchmarks motivated by subsurface geophysics,
geophysical flow forecasting, and aerodynamic design. They are not
intended to survey all classes of parametric PDEs. Rather, they span a
static elliptic problem, a time-dependent periodic flow, and a
geometry-dependent turbulent-flow problem, and are selected to explore
three potential sources of online advantage for neural operators:
hardware-efficient execution, reduction of sequential time-stepping
depth, and task reduction through direct prediction of a partial
observable. In each case, we compare the neural operator with a
baseline defined by a problem-matched classical solver. \Cref{tab:benchmark-design} 
summarizes the setup.

\begin{table}[htbp]
\centering
\small
\begin{tabular}{
@{}
>{\raggedright\arraybackslash}p{0.20\textwidth}
>{\raggedright\arraybackslash}p{0.19\textwidth}
>{\raggedright\arraybackslash}p{0.23\textwidth}
>{\raggedright\arraybackslash}p{0.27\textwidth}
@{}}
\Xhline{1.1pt}
Benchmark and application
&
Requested output and error
&
Classical baseline
&
Acceleration mechanism
\\
\hline
Darcy flow
(subsurface geophysics)
&
Pressure field; relative \(L^2\) error
&
finite elements with geometric multigrid
&
Hardware utilization  
\\ \hline
Incompressible Navier-Stokes flow
(geophysical flow forecasting)
&
Future vorticity field; relative \(L^2\) error
&
Fourier pseudospectral discretization with RK4 time integration
&
Reduction of sequential time-stepping depth
\\ \hline
Vehicle aerodynamics
(digital wind tunnel)
&
Surface pressure coefficient; relative \(L^1\) error
&
Finite-volume RANS discretization with SIMPLE iteration
&
Output reduction and avoidance of nonlinear iterations
\\
\Xhline{1.1pt}
\end{tabular}
\caption{The three benchmarks, motivating applications, requested outputs, problem-matched classical baselines, and potential sources of online advantage for neural operators.}
\label{tab:benchmark-design}
\end{table}

The Darcy benchmark tests whether neural operator inference can achieve a lower matched-accuracy cost than geometric multigrid for a coercive elliptic problem, for which the classical algorithm has nearly optimal complexity. The baseline is implemented using \texttt{Firedrake} and
\texttt{PETSc}~\cite{rathgeber2016firedrake,balay2019petsc}.
The Navier-Stokes benchmark tests whether a learned time-advance map can outperform an accurate pseudospectral solver with RK4 time integration by reducing sequential time-stepping depth. The baseline is implemented in \texttt{GeophysicalFlows.jl}~\cite{GeophysicalFlowsJOSS}.
The vehicle benchmark tests whether direct prediction of the surface
pressure coefficient can reduce online cost by avoiding both a full
volumetric solution and the nonlinear iterations of a finite-volume RANS solver. The baseline uses \texttt{OpenFOAM}
\cite{weller1998tensorial,OpenFOAM12}, a representative industrial
computational fluid dynamics (CFD) solver.

\subsection{Kernel-based neural operator}

Kernel-based neural operators approximate mappings between function spaces using layers that combine learned kernel integral operators, pointwise transformations, and nonlinear activations~\cite{li2020neural,kovachki2023neural,lanthaler2025nonlocality}.
The kernels may be localized in physical space~\cite{li2020neural},
parameterized using Fourier modes~\cite{li2020fourier},
multiwavelets~\cite{gupta2021multiwavelet}, or Laplace-Beltrami eigenbases~\cite{bonev2023spherical,chen2024learning,stuart2025enforcing}, or made input-dependent
through attention mechanisms~\cite{cao2021choose,calvello2025continuum}.
Across all three benchmarks, we use a modified Fourier neural operator (FNO)~\cite{li2020fourier} designed to accommodate problems with variable geometry~\cite{lingsch2023beyond,zeng2025point} and to better capture local effects~\cite{liu2024neural-gl,zeng2025point,han2026geometric}. 
Using the
same architecture across the benchmarks helps distinguish the influence of
problem structure from that of architecture choice. Although alternative neural operator architectures may yield different computational costs and
attainable accuracies, the acceleration mechanisms identified here are shared by many such architectures, and our conclusions are therefore expected to remain
relevant more broadly. Broader architecture comparisons are available in
\cite{lu2021comprehensive,de2022cost}, and an additional comparison with a standard FNO is provided in the supplement in~\cref{supp-sec:numerics-sensitivity}.

Given an input pair \((a,D)\), the model first lifts \(a\) to a latent field, applies \(L\) kernel integral operator layers, and then projects the final latent field to the requested output:
\begin{equation} \label{eq:NO_architecture} 
\begin{aligned} 
g_0 &= \cP(a),\\ 
g_i &= \cL_i^D g_{i-1} = g_{i-1} + \sigma\!\left( \mathcal K_{{\rm long},i}^{D}g_{i-1} + \mathcal K_{{\rm short},i}^{D}g_{i-1} \right), \qquad i=1,\ldots,L,\\ \cG(a,D;\theta) &= \cQ(g_L). 
\end{aligned} \end{equation}
Here $\cP$ and $\cQ$ are pointwise lifting and projection maps, respectively, and each latent field $g_i:D\to\R^{d_g}$ has \(d_g\) channels. Each layer \(\cL_i^D\) combines a global component \(\mathcal K_{{\rm long},i}^{D}\), which communicates information across the domain, with a local component \(\mathcal K_{{\rm short},i}^{D}\), which captures short-range spatial variation. The activation $\sigma$ is the pointwise GELU function~\cite{hendrycks2016gaussian}. The residual formulation is used as it facilitates the optimization of deeper architectures~\cite{he2016deep,WeinanE20171,haber2018stable}.

\paragraph{Global interaction}
 Suppressing the layer index, the global component is chosen to be a truncated Fourier integral operator~\cite{li2020fourier,kovachki2023neural,nelsen2021random,nelsen2024operator,kossaifi2024library,lanthaler2025nonlocality}:
\begin{equation}
\label{eq:FNO-kernel}
\begin{split}
(\mathcal K_{\rm long} g)(x)
&=
\sum_{\|k\|_\infty\le k_{\max}}
\int_{D}\kappa_{k}(x,y)W_k^v\,g(y)\,dy \\
&=
\sum_{\|k\|_\infty\le k_{\max}}
\phi_k(x)\,
W_k^v
\int_{D}\overline{\phi_k(y)}\,g(y)\,dy,
\end{split}
\end{equation}
where \(k\in\mathbb Z^d\) denotes a Fourier mode, $\kappa_k(x,y)=
e^{2\pi i k\cdot \frac{x-y}{l}}$, $\phi_k(x)=e^{2\pi i k\cdot x/l}$,
\(l\) is a characteristic length scale of the domain, and
\(W_k^v\in\mathbb C^{d_g\times d_g}\) is the learned complex matrix associated with mode
\(k\). 
Only modes satisfying \(\|k\|_\infty\le k_{\max}\) are retained. This truncation reduces the number of learnable parameters and yields a
low-rank spectral representation of the global operator.
When the computational domain $D$ is a fixed periodic box discretized on a structured grid,
\eqref{eq:FNO-kernel} can be evaluated efficiently using the fast Fourier transform (FFT):
\begin{align}
\label{eq:FNO}
    \mathcal{K}_{\rm long} g = \cF^{-1} \Bigl(R \cdot \bigl(\cF g\bigr)\Bigr),  \qquad R(k)
=
\begin{cases}
W_k^v, & \|k\|_\infty\le k_{\max},\\
0, & \text{otherwise}.
\end{cases}
\end{align}
Here \(\cF\) and \(\cF^{-1}\) denote the channelwise Fourier transform and its inverse. Thus, the global component applies learned channel mixing to each retained mode. For the Darcy flow and Navier-Stokes benchmarks, the global component
is evaluated using the FFT formulation~\eqref{eq:FNO}. For the vehicle aerodynamics benchmark, \(D\) is a variable, unstructured
surface, and the integral representation~\eqref{eq:FNO-kernel} is
evaluated directly by quadrature.

\paragraph{Local interaction}
The local component combines a pointwise affine map with a gradient correction~\cite{zeng2025point,liu2024neural-gl}:
\begin{equation}
\label{eq:smoothed-gradient}
    \mathcal{K}_{\rm short} g(x)  =  
    W g(x) + b + W^{g'}\texttt{SoftSign}\Bigl(W^{g}\textrm{vec}\bigl(\nabla_D g(x)\bigr) \Bigr).
\end{equation}
Here \(W\in\mathbb R^{d_g\times d_g}\), \(b\in\mathbb R^{d_g}\),
\(W^g\in\mathbb R^{d_g\times d_g d}\), and $W^{g'}\in\R^{d_g\times d_g}$ are learned parameters. The correction supplies local derivative information that complements the truncated global representation. The componentwise map $\texttt{SoftSign}(z)=z/(1+|z|)$ controls large discrete gradients near sharp features. We approximate \(\nabla_D g\) by centered differences on structured grids and by local least squares on point clouds; on an embedded surface, the latter approximates the tangential gradient.

\paragraph{Training} 
The lifting, projection, and layerwise parameters are collected in $\theta$. 
Given $N$ independent training inputs
\begin{equation}
    (a_j,D_j)\stackrel{\mathrm{i.i.d.}}{\sim}\mu,
    \qquad j=1,\ldots,N,
\end{equation} and corresponding reference outputs
\(\omap(a_j,D_j)\), we define the empirical measure
 $\mu^N := \frac1N\sum_{j=1}^N \delta_{(a_j, D_j)}$, and the empirical relative error loss
\begin{equation}
\label{eq:emp_risk}
    \begin{split}
                           \cJ_{N}(\theta):=&\E_{(a,D)\sim\mu^N}\Bigl( \frac{\|\omap(a,D)-\cG(a,D;\theta)\|_{\cU}}{\|\omap(a,D) \|_{\cU}}\Bigr) \\ =& \frac{1}{N}\sum_{j=1}^N\frac{\|\omap(a_j,D_j)-\cG(a_j,D_j;\theta)\|_{\cU}}{\|\omap(a_j,D_j)\|_{\cU}}. 
    \end{split}
\end{equation}
Here $\|\cdot\|_\cU$ denotes the benchmark-specific \(L^2\) or \(L^1\) norm listed in
\cref{tab:benchmark-design}, with dependence on
\(D\) suppressed. For the time-dependent Navier-Stokes benchmark, we adapt
the training loss to employ the recurrent rollout relative-error loss defined in \cref{eq:loss_rollout}.

\subsection{Measuring accuracy and cost}
We consider finite sets of classical solver and neural operator configurations and report all resulting cost-error pairs. The configurations used for each benchmark are detailed in \cref{sec:num_study}. Our accounting of online cost begins after each input has been discretized and converted to the representation required by the corresponding implementation. Accordingly, the reported costs exclude geometric preprocessing and mesh generation. For each configuration, per-query cost is quantified separately by floating-point work and wall-clock runtime.

\paragraph{Floating-point work}
We count each multiplication or addition as one flop and ignore memory traffic, communication, and parallelization overhead. For the classical methods, the count includes query-dependent assembly and iterative solution or time integration. For the neural operator, it includes lifting, all operator layers, and output projection. 
The classical solver flop counts are given with the corresponding benchmarks in \cref{sec:num_study}; the neural operator flop counts are given below.

Throughout this work, \(n_e\) denotes the effective spatial discretization
size associated with \(D\). Depending on the discretization, it counts grid
points, finite element degrees of freedom, or points in an unstructured
point cloud.
Let \(K=(2k_{\max}+1)^d\) denote the number of retained Fourier modes. For fixed input and output dimensions and retaining the leading terms, the inference count for a configuration with latent width \(d_g\) and \(L\) layers is
\begin{align}
\label{eq:no-cost}
 C_{\rm NO}^{\rm infer}(n_e,K,d_g,L)
 &\approx 8 L K d_g^2
 + \bigl(2 + (2d+4)L\bigr)d_g^2 n_e \notag\\
 &\quad+
 \begin{cases}
 10 L d_g n_e \log_2 n_e, & \text{structured grid},\\
 12 K L d_g n_e, & \text{unstructured grid}.
 \end{cases}
\end{align}
The case-dependent final term accounts for transformations between physical
and spectral representations in the global interaction: FFTs on a structured grid \eqref{eq:FNO} and separable quadrature \eqref{eq:FNO-kernel} on an unstructured point cloud. The first term accounts for modewise channel mixing, and the second collects the dominant pointwise operations and gradient corrections. For fixed \(K\), \(d_g\), and \(L\), the structured-grid cost is linear in \(n_e\), with a logarithmic correction; the unstructured-grid cost is linear in \(K\,n_e\). Detailed counts and implementation choices are provided in the supplement in \cref{supp-sec:NO}. 
Discussion of the choice of $k_{\max}$, and hence of $K$, may be found in \cite{lanthaler2025nonlocality}.

\paragraph{Wall-clock runtime}
To reflect practical computational workflows, each classical solver is
evaluated in its native software environment using the resource allocation
specified in~\cref{tab:benchmark-runtime}. All timings are obtained on compute nodes equipped with Intel Xeon Platinum 8358 CPUs and NVIDIA A100 GPUs. CPU and GPU timings are reported separately whenever both implementations are available. Unlike floating-point operation counts, Wall-clock runtime reflects memory access, parallelism, implementation, and hardware utilization in addition to floating-point work. Because resource allocations differ across methods, the reported timings should be interpreted as time-to-solution measurements for the specified implementations, hardware, and resource allocations.

\begin{table}[htbp]
\caption{Software and computational resources used for the wall-clock
evaluation of the classical solver baselines.}
\label{tab:benchmark-runtime}
\centering
\small
\begin{tabular}{
@{}
>{\raggedright\arraybackslash}p{0.20\textwidth}|
>{\raggedright\arraybackslash}p{0.19\textwidth}
>{\raggedright\arraybackslash}p{0.23\textwidth}
>{\raggedright\arraybackslash}p{0.27\textwidth}
@{}}
\Xhline{1.1pt}
Benchmark
&
Darcy flow
&
Incompressible Navier-Stokes flow
&
Vehicle aerodynamics
\\
\hline

Baseline software
&
\texttt{Firedrake/PETSc} (Python/C)
&
\texttt{GeophysicalFlows.jl} (Julia)
&
\texttt{OpenFOAM} 12 (C++)
\\
\hline

Baseline resource allocation
&
One CPU core
&
One CPU core or one GPU
&
4, 32, or 256 CPU cores, with approximately \(10^5\) cells per core
\\
\Xhline{1.1pt}
\end{tabular}
\end{table}

\section{Three cost-accuracy regimes}
\label{sec:num_study}

In this section, we use the three benchmarks to examine three regimes in which neural operators may achieve lower online cost at matched accuracy. Each benchmark is organized in the same way: we first introduce the problem-matched classical method and its expected cost-accuracy scaling, then compare the measured cost-accuracy curves and identify the mechanism underlying any observed neural operator advantage, which typically occurs for low-to-moderate accuracy demands, and finally discuss the principal limitation of neural operators revealed by the experiment, typically apparent at higher accuracy.

\subsection{Elliptic problems with optimal solvers: Darcy flow}
Darcy flow is a stringent test for neural operator inference because geometric multigrid has nearly linear complexity for the resulting coercive elliptic problem. The comparison therefore distinguishes a reduction in floating-point work from an advantage in hardware execution.

\begin{figure}[htbp]
     \centering
     \includegraphics[width=0.99\textwidth]{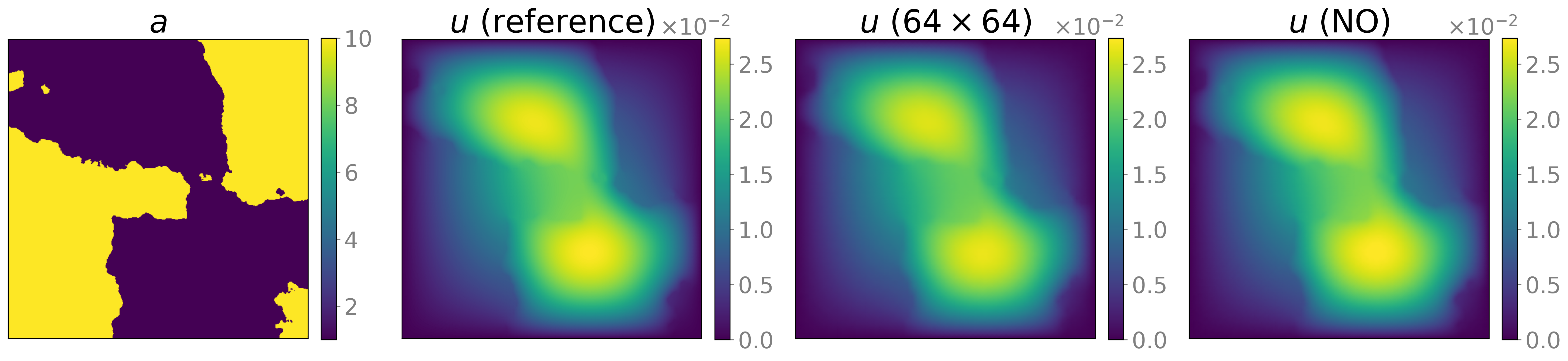}\\
     \caption{Darcy flow for a representative test sample.
     From left to right: the permeability field, the reference solution, the FEM solution computed with \texttt{Firedrake} on a coarse \(64 \times 64\) grid, and the neural operator (NO) prediction.}
     \label{fig:darcy-map}
\end{figure}

We consider the two-dimensional Darcy flow problem on the square domain $D=[0,1]^2$:
\begin{equation}
\begin{aligned}
    -\nabla \cdot (a \nabla u) &= f \qquad &&\text{in } D, \\
    u &= 0 \qquad &&\text{on } \partial D .
\end{aligned}
\label{eq:Elliptic}
\end{equation}
The forcing term is fixed as $f\equiv 1$. Our goal is to learn the solution operator 
\begin{equation}
    \omap : \bigl(a, [0,1]^2\bigr) \mapsto u;
\end{equation}
the domain $D$ is fixed across all samples, is encoded in the neural operator through the coordinate map $x \mapsto x$ and we do not attempt to learn dependence on varying $D$.
The permeability field, illustrated in \cref{fig:darcy-map}, is discontinuous and is modeled as
\[
a(x) = 1 + 9 \mathbf{1}_{\{\tilde a(x) \ge 0\}},
\]
where $\tilde a$ is a centered Gaussian random field with covariance operator
\[
\mathcal C = (-\Delta + \tau^2)^{-\alpha}.
\]
Here $-\Delta$ denotes the Laplacian on $D$, subject to Neumann boundary conditions on the space
of spatial-mean zero functions, with $\alpha=2$ and $\tau=3$ following \cite{li2020fourier}. Thus, $a$ takes values in $\{1,10\}$, giving a higher contrast permeability field than in the original setup of \cite{li2020fourier}.

\subsubsection{Classical baseline and expected scaling}
We discretize the Darcy flow equation~\eqref{eq:Elliptic} using conforming bilinear \(Q_1\) finite elements on a uniform
\(n \times n\) quadrilateral mesh of \(D=[0,1]^2\), with mesh size \(h=1/n\) and $n_e=n^2$ elements.
Let \(V_h\) denote the corresponding finite element space with homogeneous Dirichlet
boundary conditions. The discrete solution \(u_h \in V_h\) satisfies
\[
\int_D a \nabla u_h \cdot \nabla v_h \, dx
=
\int_D f v_h \, dx,
\qquad \forall v_h \in V_h .
\]
Both the permeability field \(a\) and the source term \(f\) are represented by nodal
\(Q_1\) functions on the same mesh, and the integrals are evaluated using tensor-product
Gaussian quadrature. The resulting discretization yields a sparse symmetric positive
definite linear system
\[
A_h u_h = b_h.
\]
We solve this system using geometric multigrid as
a standalone iterative solver, specifically a multiplicative \(V\)-cycle
\cite{briggs2000multigrid,trottenberg2001multigrid} on a hierarchy of uniformly refined quadrilateral
meshes. On each level, we use Richardson-Jacobi smoothing, while a direct solver is
applied on the coarsest grid. 
For the implementation, we use \texttt{Firedrake} \cite{rathgeber2016firedrake} as a robust and reproducible finite element baseline, with \texttt{PETSc} \cite{balay2019petsc} providing the underlying linear algebra and multigrid infrastructure. 
The reference data are generated on a \(512 \times 512\) grid. The linear iterations are continued until the residual reaches machine precision.

The cost of the finite element multigrid solver can be decomposed into an assembly
cost and an iterative solve cost. The assembly phase includes the construction of the stiffness matrices, the finest-grid load vector, and the intergrid transfer operators. 
The solve phase consists of \(m\) multigrid
\(V\)-cycles. Up to lower-order terms, the
resulting floating-point cost satisfies
\begin{equation}
\label{eq:darcy_C_total}
C_{\rm solver}
\approx
\left(
c_{\rm assemble}
+
c_{\rm solve}\, m
\right) n_e.
\end{equation}
The constants \(c_{\rm assemble}\) and \(c_{\rm solve}\) are both on the order of a few hundred; their detailed values are reported in the supplement in \cref{supp-ssec:darcy}.
If the coarsest grid is fixed and the number of \(V\)-cycles required to reach a fixed
relative residual tolerance is independent of the mesh size, so that $m = \bigO(1)$, 
then the total solver cost is
linear in the number of finest-grid degrees of freedom:
\[
C_{\rm solver} = \mathcal O(n_e).
\]
Mesh-independent convergence of multigrid under standard smoothing and approximation
assumptions is classical; see, for example,
\cite{braess1983new,xu1992iterative,xu2002method}.
For sufficiently regular solutions, the \(Q_1\) finite element method satisfies the error estimate
\[
\|u-u_h\|_{L^2}=\mathcal O(h^q),
\qquad q \le 2.
\]
In the smooth coefficient case, the optimal rate $q=2$ is recovered. In the present benchmark, however, the permeability field is discontinuous, and the effective convergence rate may therefore be lower~\cite{babuvska2000can}. Assuming that the algebraic error is reduced below
the discretization error, achieving an \(L^2\) error of order \(\varepsilon\) requires \(h=\bigO(\varepsilon^{1/q})\). Since \(n_e=\bigO(h^{-2})\) in two dimensions, the corresponding cost scales as
\begin{equation}
\label{eq:darcy_cost_accuracy_theory}
    C_{\rm solver}(\varepsilon) = \mathcal O(n_e) = \mathcal O(h^{-2}) = \mathcal O(\varepsilon^{-2/q}).
\end{equation}

\subsubsection{Measured trade-off and acceleration mechanism}
For the classical numerical solver, we use the grid sequence $n_e = n\times n, \, n = 16,\ 32, \ldots,\ 512$, with an inexact solve terminated at a relative residual tolerance of \(10^{-6}\). 
The neural operator uses the structured-grid implementation with the zero padding strategy for the treatment of nonperiodic boundaries described in the supplement in \cref{supp-rem:bc}. The input has $d_a = 3$ channels: the permeability field $a(x)$ and two coordinate functions $x_1$ and $x_2$.
We fix the training set size at \(N=4000\), determine the number of
Fourier modes by setting \(k_{\max}=16\) and set the latent width to 
\(d_g=64\). We vary the downsampled grid resolution over $n = 32,\ 64$, and  $128,$
and the number of neural operator layers over $L = 4,\ 5,\ 6 $, to obtain a range of
 surrogate implementations.

\begin{figure}[htbp]
     \centering
     \includegraphics[width=0.9\textwidth]{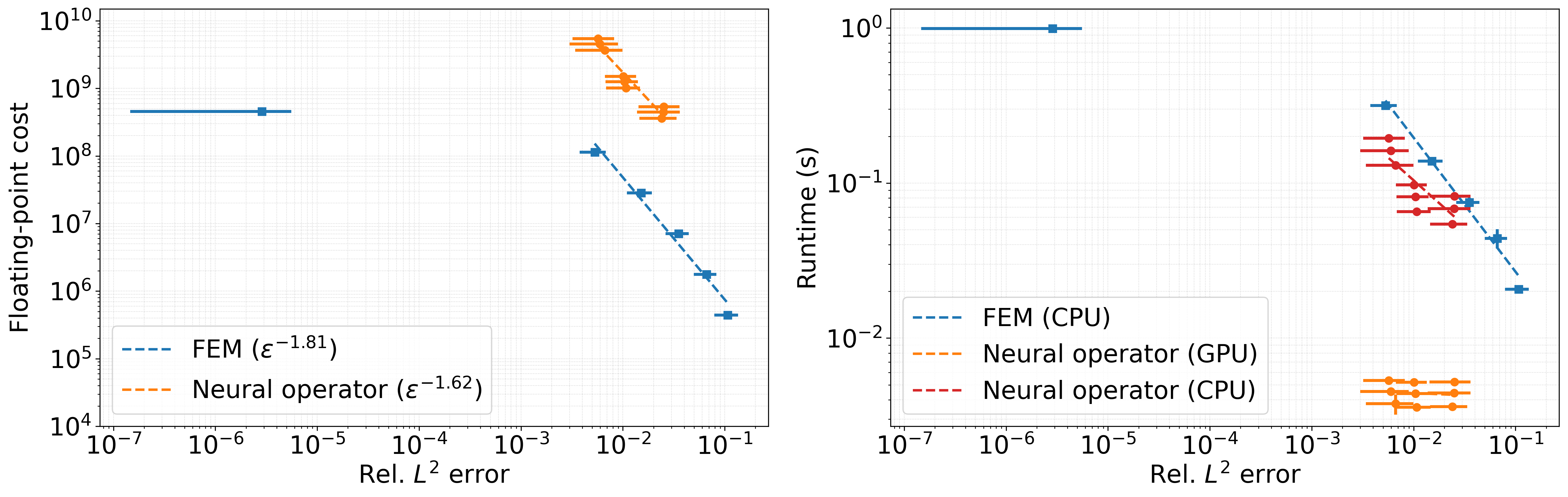}\\
     \caption{Cost-accuracy comparison for Darcy flow between a finite element method with geometric multigrid (FEM) and the neural operator.
     Left: floating-point operation count versus relative $L^2$ error. 
     Right: wall-clock runtime in seconds versus relative $L^2$ error for CPU and GPU 
     implementations. Error bars denote one standard deviation over the test samples.
     The isolated FEM point at the smallest reported error is obtained on the same \(n=512\) grid as the reference solution but with fewer solver iterations. Because this comparison measures algebraic iteration error but not spatial discretization error, the point is excluded from the FEM power law fit.}
     \label{fig:darcy-cost-accuracy}
\end{figure}

The cost-accuracy comparison is reported in \cref{fig:darcy-cost-accuracy}. 
In terms of floating-point operation count, the FEM with geometric multigrid is highly efficient in the high accuracy regime: it typically converges in fewer than 20 multigrid iterations and exhibits nearly optimal scaling, with cost
linear in \(n_e\), as shown in \cref{eq:darcy_C_total}. 
Its empirical cost-accuracy curve corresponds to $q \approx 1.1$ in \cref{eq:darcy_cost_accuracy_theory}, which is similar to that of the neural operator
with $q \approx 1.2$ (leading to the predicted growth of cost as a function of $\varepsilon$
in the figure). However the cost-accuracy curve for the FEM method lies far below that of
the neural operator, and is far less noisy.
The wall-clock comparison leads to a different conclusion:
at comparable low-to-moderate accuracy, the neural operator is slightly faster on CPUs and achieves an additional order-of-magnitude speedup on GPUs. This acceleration is largely due to the highly modular structure of the neural operator, whose inference consists primarily of matrix-vector products, FFTs, and pointwise nonlinear activations. These operations can be executed efficiently on modern accelerator hardware. Thus, this benchmark exhibits an implementation- and hardware-execution advantage rather than a reduction in floating-point work.


\subsubsection{Limitations and implications}
The neural operators in \cref{fig:darcy-cost-accuracy} attain relative errors between \(5\times10^{-3}\) and \(10^{-1}\). To examine whether additional data
or model capacity can further reduce this error, we conduct two studies on data downsampled on a \(128\times128\) grid.

First, we fix \(k_{\max}=16\), \(L=6\), and \(d_g=64\), and vary the training-set size \(N\). The corresponding errors are reported in \cref{tab:darcy}-top. Error decreases systematically with \(N\), but the improvement slows as the dataset grows. A log-log fit to the displayed values
gives an empirical decay of approximately \(N^{-0.30}\). Second, we fix \(N=8000\) and compare architectures with different spectral
cutoffs $k_{\max}$, widths $d_g$, and depths $L$, as shown in \cref{tab:darcy}-bottom. Increasing model capacity improves accuracy,
but the gains diminish once the error approaches \(5\times10^{-3}\). Published results on related Darcy benchmarks similarly report errors of order \(10^{-3}\) or larger for standalone neural operators~\cite{li2020fourier,he2024mgno}. These
results indicate an empirical error plateau for the tested configurations and
training procedure. Therefore, the tested neural operators do not reach the higher accuracy
regime accessible to the finite element solver through mesh refinement.

\begin{table}[htbp]
    \begin{center}
        \begin{tabular}{c|cccc}
        \Xhline{1.1pt}
            $N$  & $1000$ & $2000$ & $4000$ & $8000$\\ \hline
            Rel. $L^2$ test error~($\times 10^{-3}$)  & $9.18$   &   $6.82$ & $5.56$ & $4.91$ \\\Xhline{1.1pt}
            ($k_{\max}, d_g, L$)   & ($4, 16, 4$) & ($8, 32, 5$) & ($16, 64, 6$)  & ($32, 128, 7$)\\ \hline
            Rel. $L^2$ test error~($\times 10^{-3}$) & $11.37$  & $5.46 $   &   $4.91$ & $4.70$ \\\Xhline{1.1pt}
        \end{tabular}
    \end{center}
    \caption{Darcy flow: test error versus training set size $N$ for a fixed neural operator architecture \(k_{\max}=16\), \(L=6\), \(d_g=64\) (top); test error versus architecture ($k_{\max}, d_g, L$) for a fixed training data size $N=8000$ (bottom). }
    \label{tab:darcy}
\end{table}

A useful heuristic decomposition of the in-distribution test error for a neural operator is
\cite{lu2021learning,de2022cost}
\begin{equation}
\label{eq:empirical_predictive_error}
\E_{(a,D)\sim\mu}\Bigl( \frac{\|\omap(a,D)-\cG(a,D;\theta_N)\|_{\cU}}{\|\omap(a,D) \|_{\cU}}\Bigr) 
\approx
c_{\rm disc} h_{\rm eff}^q + c_{\rm data} N^{-\beta} + \varepsilon_{\rm opt}.
\end{equation}
Here \(\theta_N\) denotes neural operator parameters trained on \(N\) samples.
The first term represents discretization and architectural representation errors. The effective scale \(h_{\rm eff}\) may represent the physical grid spacing
when spatial resolution is limiting, or \(k_{\max}^{-1}\) when Fourier
truncation is dominant. Encoder-decoder truncation produces analogous
representation errors. Approximation results that motivate these contributions
include \cite[Theorem~24]{kovachki2021universal}
\cite[Theorem~2.1]{han2026geometric}
\cite[Theorem~3.4]{bhattacharya2021model}\cite{schwab2026deep}. The exponent \(q>0\) is therefore an
effective approximation exponent that depends on both the regularity of the target operator and the approximation properties of the architecture.
The second term represents the statistical error arising from the approximation of the population loss
\begin{equation}
\label{eq:pop_risk}
    \begin{split}
                           \cJ_{\infty} (\theta):=&\E_{(a,D)\sim\mu}\Bigl( \frac{\|\omap(a,D)-\cG(a,D;\theta)\|_{\cU}}{\|\omap(a,D) \|_{\cU}}\Bigr)  
    \end{split}
\end{equation}
by the empirical loss~\eqref{eq:emp_risk} based on \(N\) training samples. The exponent \(\beta>0\) generally depends on the regularity of the target solution operator and on the input distribution; see, for example, \cite[Theorem 1.3]{de2023convergence}\cite[Corollary 2.8]{reinhardt2026statistical}\cite[Theorem 4.9]{huang2024operator}\cite{liu2024neural}\cite[Theorem 4]{liu2024deep}. 
The final term, \(\varepsilon_{\rm opt}\), accounts for residual error because
training only approximately minimizes the nonconvex empirical objective. This
residual may arise from finite optimization time, stochastic gradient noise, or
convergence to a suboptimal stationary point
\cite{bottou2007tradeoffs,cisneros2025optimization}.

The data and architecture studies show diminishing improvements, suggesting that
fixed resolution representation error or residual optimization error may now be
limiting.
Achieving substantially smaller errors with a stand-alone neural
operator therefore becomes increasingly difficult,  although recent studies suggest
that second-order optimizers can improve attainable accuracy in related
scientific machine learning models~\cite{kiyani2025optimizing,
jnini2026curvature}. For applications requiring
stringent accuracy, hybrid methods may offer a more promising path: the neural
operator supplies an initial guess, coarse correction, or component of a
preconditioner, while a classical solver delivers the final high accuracy
solution~\cite{zhang2024blending,hu2025hybrid}.

\subsection{Time-dependent problems: incompressible Navier-Stokes flow}

\begin{figure}[htbp]
     \centering
     \includegraphics[width=0.9\textwidth]{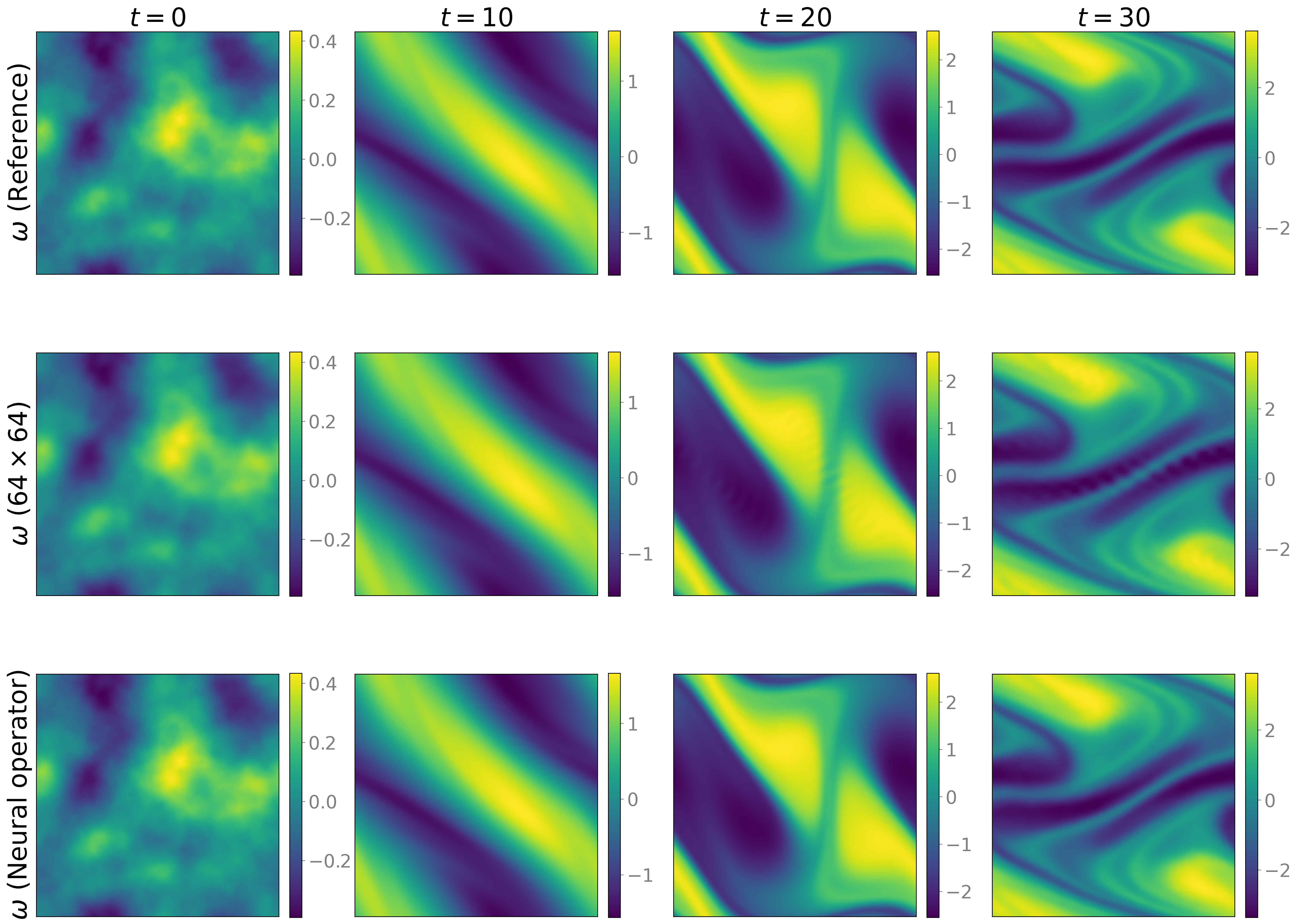}\\
     \caption{Incompressible Navier-Stokes flow for a representative test trajectory at \(t=0,10,20,\) and \(30\) (columns). From top to bottom: the \(256\times256\) reference solution, the spectral solution computed with \texttt{GeophysicalFlows.jl} on a coarse $64\times64$ grid, and the neural operator rollout prediction. }
     \label{fig:NS-map}
\end{figure}

This time-dependent benchmark tests a different mechanism. To advance the solution by one unit of physical time, the
explicit pseudospectral solver performs many CFL-limited RK4 steps, whereas the
neural operator applies one learned time-advance map. The potential advantage
is reduced sequential depth.
 
We consider the two-dimensional incompressible Navier-Stokes equations on the periodic
square domain \([0,1]^2\); thus, $D=\mathbb{T}^2$, the two-dimensional unit torus. In vorticity form, the equation is
\[
\partial_t \omega + v\cdot \nabla \omega
=
\nu \Delta \omega + f,
\]
where $\omega$ is the scalar vorticity, $\nu$ is the viscosity, and \(v\)
is the divergence-free velocity recovered from \(\omega\) through the periodic
Biot-Savart law. We consider fixed vorticity forcing  $f_{\mathsf{bench}}(x) = 0.1\bigl(\sin(2\pi (x_1 + x_2)) + \cos(2\pi (x_1 + x_2))\bigr)$, following \cite{li2020fourier}. 
Under this Kolmogorov-type forcing, the flow becomes chaotic after an initial transient; see \cref{fig:NS-map}.
The goal is to learn the unit-time solution operator
\[
\omap  : \bigl(\omega(t), f_{\mathsf{bench}}, \mathbb{T}^2\bigr) \mapsto \omega(t+1).
\]
By iterating the learned map, we obtain a recurrent surrogate for predicting flow transitions and certain chaotic features over longer time horizons. 
In this problem,
both the vorticity forcing \(f\) and the domain \(D\) are fixed across all samples, but they are included as input features to enrich the representation.
The initial vorticity $\omega_0$ is modeled as a centered Gaussian random field with covariance
operator
\[
\mathcal C = \tau^{2\alpha - d}(-\Delta+\tau^2)^{-\alpha},
\]
where \(-\Delta\) denotes the Laplacian on periodic, spatial-mean-zero functions. Following \cite{li2020fourier}, we take
\(\alpha=2.5\), \(\tau=7\), $\nu=10^{-4}$. With this choice, \(\omega_0\in H^{q_0}(D)\) almost surely for every \(q_0<\alpha-d/2=1.5\).

\subsubsection{Classical baseline and expected scaling} 

We employ a Fourier pseudospectral discretization in space~\cite{hesthaven2007spectral}, together with a fourth-order
Runge-Kutta method (RK4) in time. The domain is discretized on a uniform periodic \(n \times n\) grid with \(n_e = n^2\) points and mesh size \(h=1/n\). In this setting, the Fourier coefficients \(\widehat{\omega}_k\) of the vorticity satisfy
\[
\partial_t \widehat{\omega}_k
=
-\widehat{(v\cdot\nabla\omega)}_k
-\nu |2\pi k|^2 \widehat{\omega}_k
+\widehat f_k,
\]
where \(k\) denotes the Fourier index. 
Spatial derivatives, the Laplacian, and the stream-function Poisson solve used to recover the velocity through the Biot-Savart relation, are diagonal in Fourier space and can therefore be evaluated efficiently using FFTs, whereas the nonlinear advection term is evaluated pseudospectrally.
The time step is chosen to satisfy both the advective CFL condition and the explicit viscous stability condition. Compared with lower-order explicit schemes, RK4 typically permits a larger time step for a prescribed accuracy.
A comparison with exponential time-differencing schemes~\cite{cox2002exponential} is provided in the supplement in \cref{supp-sec:numerics-etd}; the qualitative cost-accuracy conclusions remain unchanged.

For the implementation, we use \texttt{GeophysicalFlows.jl}~\cite{GeophysicalFlowsJOSS}, a Julia package for Fourier-based pseudospectral solvers on periodic domains. The package supports both CPU and GPU execution. The reference data are generated on a \(256\times 256\) grid, with a sufficiently small time step \(\Delta t=10^{-3}\). For each sampled initial vorticity field \(\omega_0\), we simulate the trajectory up to \(T=50\) and store snapshots at integer times.

The cost of the explicit time-dependent solver is the product of the number of time
steps \(n_t\) and the cost per step. The latter is determined by the four
RK4 stages and the cost of evaluating the right-hand side. For the present
pseudospectral method, the dominant cost at each stage arises from the FFT-based evaluation of
the right-hand side. Up to lower-order terms, the total floating-point operation count satisfies
\begin{equation}
    \label{eq:ns_c_total}
    C_{\rm solver}
\approx
n_t\left(c_{\rm assemble}\, n_e\log_2 n_e + c'_{\rm assemble}\, n_e\right),
\end{equation}
where \(c_{\rm assemble}\) and \(c'_{\rm assemble}\) denote the leading FFT and non-FFT contributions, respectively. Both constants are on the order of a few hundred, and their detailed values are given in the supplement in \cref{supp-ssec:ns}.
At a fixed output time \(T\), the fully discrete error contains spatial and temporal components:
\[
\|\omega(T)-\omega_{n,\Delta t}(T)\|_{L^2}
=
\mathcal O\!\left(h^{q_T} + \Delta t^4\right),
\]
where $\omega_{n,\Delta t}(T)$ is computed with mesh size $h$ and time step $\Delta t$, and \(q_T\) denotes the effective spatial convergence
exponent over the resolutions considered.
The spatial error reflects both the
Fourier pseudospectral discretization and the regularity of the solution. Although the initial vorticity has only finite Sobolev regularity, viscous smoothing increases its regularity at positive times and may yield spectral
spatial convergence~\cite{foias1989gevrey}.
The RK4 method contributes the 4th-order temporal error $\mathcal O(\Delta t^4)$.
In the present setting, the viscosity is small, and the advective CFL condition dominates, so for fixed final time $T$, 
\[
\Delta t = \bigO(h)=\bigO(n^{-1}), \qquad n_t = \bigO(n)= \bigO(n_e^{1/2}).
\]
If the observed error behaves as
\(\varepsilon = \bigO\bigl(n_e^{-\frac{\min\{q_T, 4\}}{2}}\bigr)\), the corresponding 
cost scales as
\begin{equation}
\label{eq:ns_cost_accuracy_theory}
C_{\rm solver}(\varepsilon)
=
\mathcal O\!\left(n_e^{3/2}\log_2 n_e\right)
=
\mathcal O\!\left(\varepsilon^{-3/\min(q_T,4)}\log_2 \varepsilon^{-1}\right).
\end{equation}

\subsubsection{Measured trade-off and acceleration mechanism}
For the classical numerical solver, we use the grid sequence $n_e = n\times n, \, n = 32,\ 64, \ldots,\ 256$, with time step $\Delta t = 1/(2n)$.
The neural operator uses the structured-grid implementation, adapted to periodic boundary conditions. The input has $d_a = 4$ channels: the current vorticity field \(\omega(x)\), the vorticity forcing \(f(x)\), and two coordinate functions \(x_1\) and \(x_2\).  We train on 200 trajectories up to \(T=50\), sampled at integer times and yielding approximately \(10{,}000\) training data pairs. We fix the rollout horizon~\eqref{eq:loss_rollout} at \(s=2\), \(k_{\max}=16\), and \(d_g=64\), and vary the downsampled grid resolution over $n = 32$, $64$, and $128$, the number of neural operator layers over \(L=4,5,6\), as for the Darcy example.

\begin{figure}[htbp]
     \centering
     \includegraphics[width=0.9\textwidth]{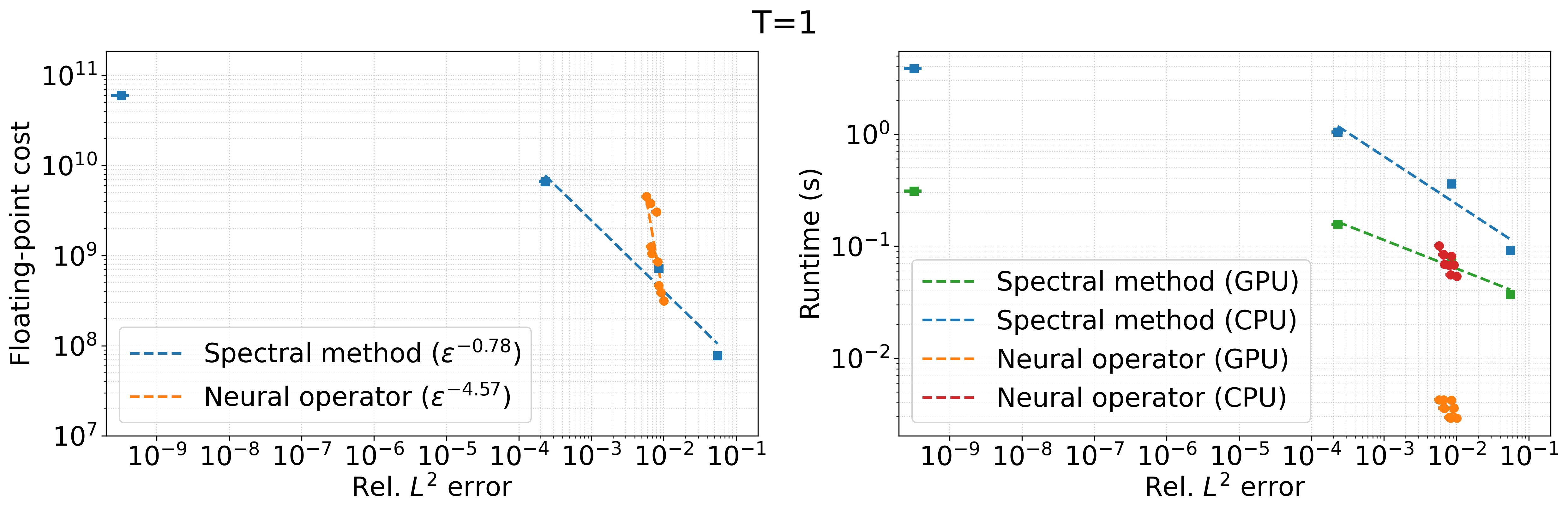}\\
     \includegraphics[width=0.9\textwidth]{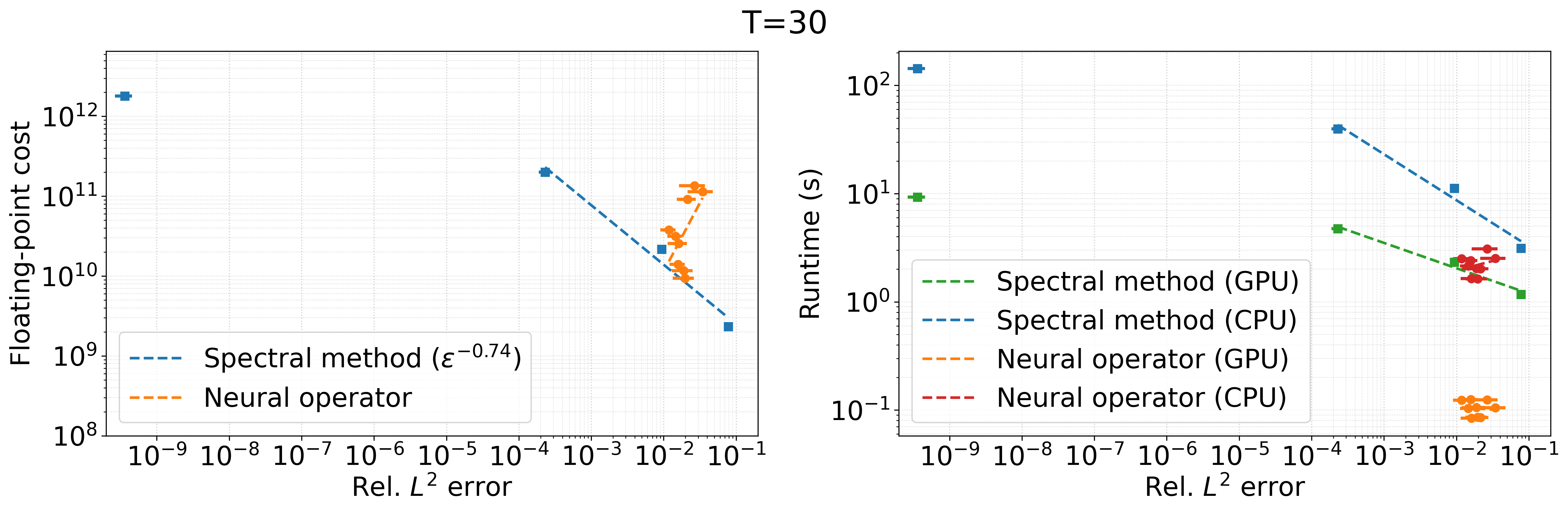}\\
     \caption{Cost-accuracy comparison for incompressible Navier-Stokes flow between
            the Fourier pseudospectral solver and the neural operator.
            Top: teacher-forced one-step evaluation \((T=1)\), with each reference
    state \(\omega(t)\), \(t=0,\ldots,29\), used to
    predict \(\omega(t+1)\). Bottom: recurrent rollout from \(t=0\)
    through \(T=30\).
     Left: floating-point operation count versus relative $L^2$ error. 
     Right: wall-clock runtime in seconds versus relative $L^2$ error for CPU and GPU 
     implementations. Error bars denote one standard deviation over the test samples.
    The isolated smallest error result for the Fourier pseudospectral solver uses
the reference grid (\(n=256\)) with a larger time step. 
Because it primarily measures temporal error rather than the fully
discrete error, it is excluded from the power-law fit.}
     \label{fig:NS-cost-accuracy}
\end{figure}

Following \cite{li2020fourier}, we consider two evaluation settings. For teacher-forced one-step evaluation, denoted by $T=1$, each reference state $\omega(t), t=0,\cdots 29$, is supplied as input to predict $\omega(t+1)$. For recurrent prediction, denoted by $T=30$, only $\omega(0)$ is supplied, and the learned unit-time map is applied recursively through \(T=30\). For each test trajectory, we report the discrete trajectory-averaged relative \(L^2\) error, defined as the average of the
relative errors over the 30 predicted states.
The present setting is more challenging than that in
\cite{li2020fourier}, where the history on \([0,10]\) is provided and only
\(t\in(10,30]\) is predicted.
During training, the unit-time update map is unrolled for $s=2$ steps and optimized using the
multistep recurrent loss in \eqref{eq:loss_rollout}.

The cost-accuracy comparison is reported in \cref{fig:NS-cost-accuracy}. 
In terms of floating-point operation count, the Fourier pseudospectral solver with RK4 is highly efficient. Its empirical cost-accuracy behavior is consistent with \cref{eq:ns_cost_accuracy_theory}: the fitted rates are $q_1 \approx 3.9$ for the teacher-forced one-step prediction and $q_{1:30} \approx 4.0$ for the trajectory-averaged recurrent rollout. Both
values are close to the fourth-order temporal-accuracy ceiling predicted by
the theory over the resolutions tested.
In both evaluation settings, the neural operator replaces the
CFL-restricted RK4 steps over each unit-time interval with a single learned
update. However, this reduction in sequential
depth does not yield a clear advantage in floating-point work at comparable
error levels, except possibly near \(\varepsilon=2\times10^{-2}\). Although each pseudospectral stage and each learned
update has quasilinear spatial cost, the spectral solver additionally requires
\(O(n_e^{1/2})\) time steps per unit-time interval, whereas the neural operator
has a large per-update prefactor.
The measured wall-clock runtime nevertheless remains favorable to neural inference on both the CPU and GPU implementations, with speedups of approximately one order of magnitude. Although both methods rely on FFTs and pointwise operations, the classical implementation advances sequential time steps that cannot be parallelized across time, while the neural operator uses one learned
update per unit time. The principal runtime mechanism is therefore reduced sequential depth.

\subsubsection{Limitations and implications}
\begin{figure}[htbp]
     \centering
     \includegraphics[width=0.6\textwidth]{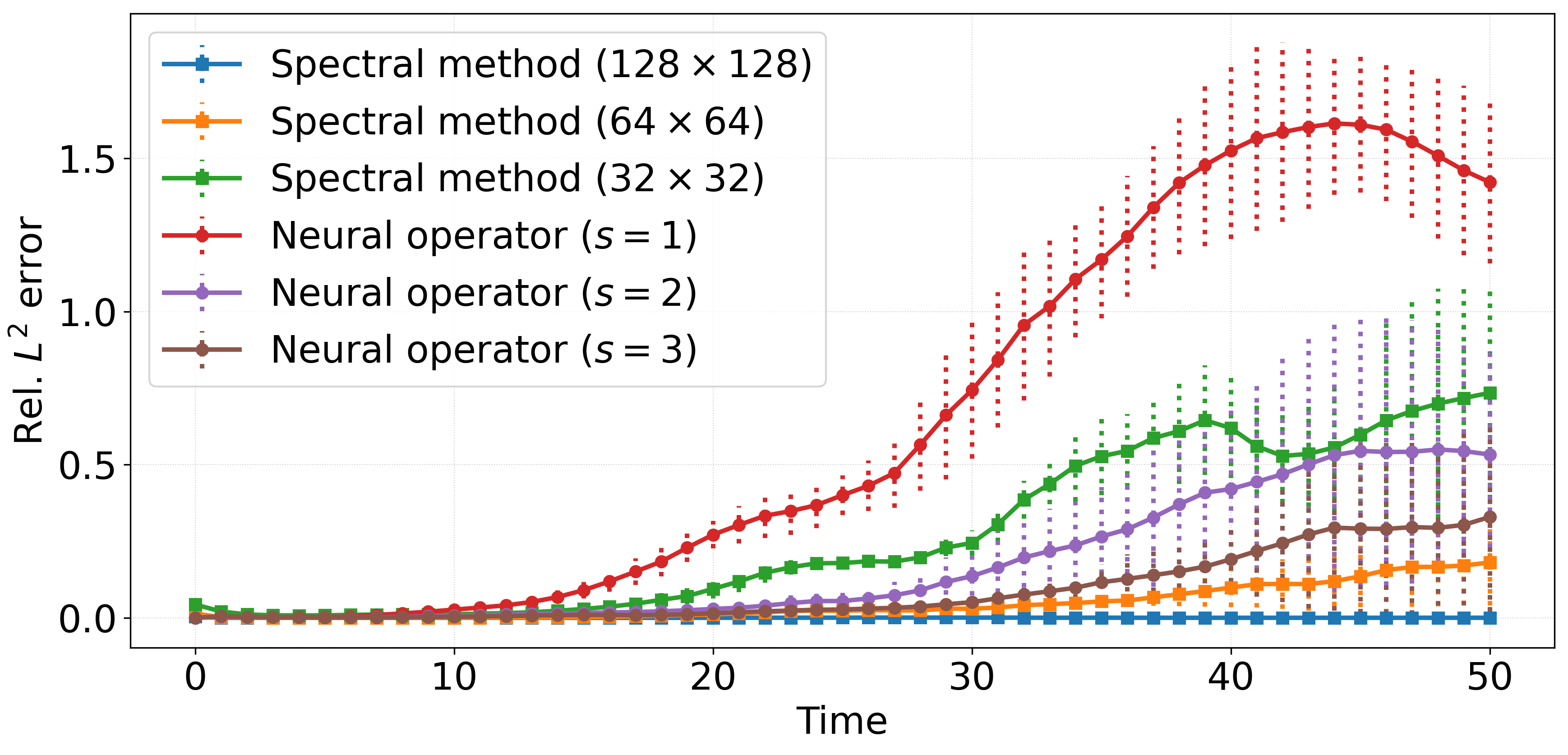}\\
     \caption{Error accumulation over time for incompressible Navier-Stokes flow,
comparing the Fourier pseudospectral solver at different spatial
resolutions with neural operator models trained using different
training horizons. Error bars denote one standard deviation over the test samples.}
     \label{fig:NS-error-accumulation}
\end{figure}

Low one-step error does not guarantee an accurate long-horizon rollout. 
In the 30-step experiment~\cref{fig:NS-cost-accuracy}, the neural operators attain an average relative \(L^2\) error of order \(10^{-2}\), but errors accumulate under recurrent application.
This growth can be mitigated by using an incremental formulation and training the update map recurrently over \(s\) steps with the multistep recurrent loss
\cite{zhou2025improving}:
\begin{align}
\label{eq:loss_rollout}
\cJ_s(\theta)
:=\E_{\substack{\omega_0\sim\mu,\\
r\sim\operatorname{Unif}\{0,\ldots,50-s\}}}
\left[
\frac1s\sum_{i=1}^{s}
\frac{\|\omega_{r+i}-\cG^{(i)}(\omega_r,f,D;\theta)\|_{L^2}}
{\|\omega_{r+i}\|_{L^2}}
\right],
\end{align}
where \(\cG^{(i)}\) denotes the \(i\)-step recurrent application of the learned
operator. Increasing \(s\) generally improves rollout accuracy, but also increases the training cost and memory requirements.
To demonstrate the effect of recurrent training, we fix the architecture at \(k_{\max}=16\), $L=6$, and $d_g = 64$, and train on data downsampled to a \(128\times 128\) grid. We compare training horizons $s=1,\ 2,\ 3$, and include the spectral solver at several grid resolutions as a reference.  
As shown in~\cref{fig:NS-error-accumulation}, the rollout error initially grows approximately linearly for all three horizons. 
Increasing \(s\) improves accuracy but does not eliminate accumulation; one possible explanation for later rapid growth is that recurrently predicted states drift away from those represented in training. In this regime, one step error is no longer a reliable indicator of rollout accuracy. 
In contrast, over the finite horizons considered here, the error of the classical solver can be reduced systematically by refining the time step and spatial grid. Although numerical errors also accumulate over time, refinement provides a controlled route to greater accuracy. For neural operators, long horizon recurrent prediction remains challenging, because errors accumulate, especially when transient dynamics must be resolved accurately. Promising directions for future work include structure-preserving
objectives based on semigroup constraints~\cite{li2022learning,chen2025due}, symplecticity-preserving architectures for Hamiltonian
systems~\cite{jin2020sympnets,makara2026symplectic},
data assimilation strategies~\cite{revach2022kalmannet,huang2025learning,haywood2026piggo},
and direct spatiotemporal prediction in latent
spaces~\cite{wiewel2019latent,hu2025deepomamba}.

\subsection{Partial-observable surrogates: vehicle aerodynamics}
\label{ssec:3D-NS}
The vehicle benchmark combines two acceleration mechanisms. The neural operator is a task-specific surrogate intended for applications such as
aerodynamic design. It predicts only the surface pressure coefficient
and does not provide the volumetric flow field, thereby bypassing both the
computation of the full volumetric flow and the nonlinear iterations required
by the classical RANS solver. 

We consider steady external aerodynamic flow around a vehicle. Let
\(B\subset\mathbb R^3\) denote the computational bounding box and let \(D\)
denote the vehicle surface. The corresponding fluid domain
\(\Omega_D\) is the portion of \(B\) exterior to the vehicle.
The flow is modeled by the steady incompressible Reynolds-averaged
Navier-Stokes (RANS) equations
\begin{subequations}
\begin{align}
    \nabla\cdot(v\otimes v)
    &=
    -\frac{1}{\rho}\nabla p
    +
    \nabla\cdot
    \left[
        (\nu+\nu_{\rm turb})
        \left(\nabla v+\nabla v^T\right)
    \right],\\
    \nabla\cdot v &= 0, 
\end{align}
\end{subequations}
where \(v\) is the mean velocity, \(p\) is the pressure, \(\rho\) is the
constant density, \(\nu\) is the molecular kinematic viscosity, and
\(\nu_{\rm turb}\) is the turbulent eddy viscosity provided by the turbulence model, here taken to be the two-equation \(k\)-\(\omega\) SST model~\cite{menter1994two,menter2003ten}.
The vehicle boundary is decomposed into stationary body surfaces, where a no-slip
condition is imposed, and wheel surfaces, where rotating wall boundary conditions are
applied. On the inlet boundary of \(B\), we prescribe a uniform incoming velocity
\(v_\infty = 30\ m\ s^{-1}\). The outlet, side, and top boundaries are treated using freestream
far-field conditions with reference velocity \(v_\infty\) and reference pressure
\(p_\infty\). The ground plane is modeled as a moving wall with velocity
\(v_\infty\), consistent with the vehicle-fixed frame of reference.

Our goal is to learn a surrogate map from the vehicle geometry 
to the surface pressure coefficient on $D$:
\begin{equation}
 \omap : D \mapsto C_p|_{D},
\quad \textrm{where} \quad 
C_p = \frac{p-p_\infty}{\frac12\rho v_\infty^2}.
\end{equation}
The boundary conditions and flow parameters are fixed across samples, so the geometry is the only varying input. This remains a parametric PDE problem, since the computational domain \(\Omega_D\) varies with the vehicle geometry.
This is a partial-surrogate setting: the neural operator predicts only the surface
quantity of interest $C_p|_{D}$, rather than the full three-dimensional flow
field or quantities like lift or drag. 
This choice is motivated by aerodynamic design, where the surface pressure distribution directly determines force and loading information and captures local geometric effects on the flow.

\begin{figure}[htbp]
     \centering
     \includegraphics[width=0.98\textwidth]{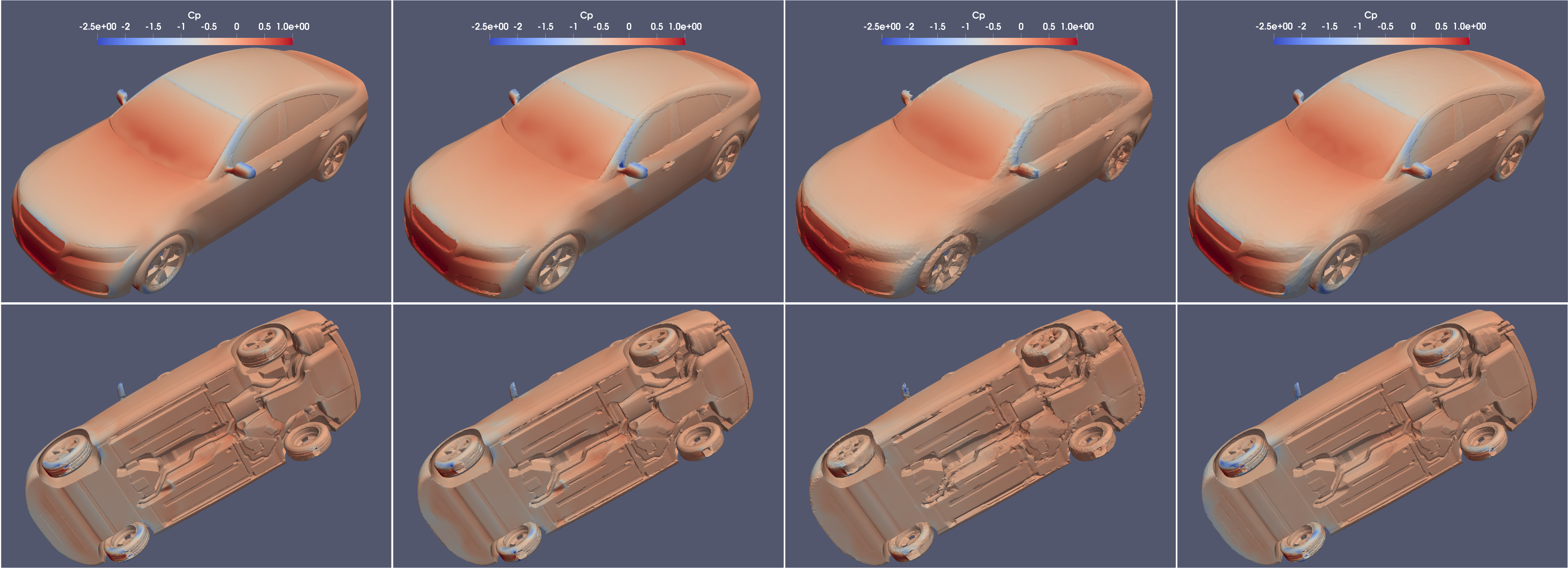}\\
     \caption{
     Surface pressure coefficient \(C_p\) for a representative vehicle from the test set.
     From left to right: the reference solution
     obtained using the large mesh and 7000 SIMPLE iterations; \texttt{OpenFOAM} predictions obtained using the medium and small meshes; and the neural operator prediction evaluated on a downsampled surface mesh containing approximately \(20{,}000\) points.
     }
     \label{fig:car-map}
\end{figure}

For variable geometries, we use the \texttt{DrivAerNet++} dataset~\cite{elrefaie2024drivaernet++},
which contains $8,000$ vehicle designs with varying geometric parameters, including
overall vehicle length and roof height, and spanning configurations such as fastback,
notchback, and estateback. The geometries exhibit complex aerodynamic features representative of industrial vehicle design; see \cref{fig:car-map}. 

\subsubsection{Classical baseline and expected scaling}

We employ a finite volume discretization in space and solve the steady nonlinear RANS
system using a segregated SIMPLE fixed-point iteration~\cite{patankar1983calculation}.
The fluid domain is discretized by a cell-centered finite volume method on an
unstructured polyhedral mesh, dominated by hexahedral cells, with \(n_c\) control
volumes. At each outer SIMPLE iteration, the momentum equations are assembled using the
current pressure, face fluxes, and turbulent viscosity, and are solved for a predicted
velocity using an iterative linear solver with symmetric Gauss-Seidel smoothing. A
pressure correction equation is then derived from the discrete continuity constraint and
solved using a geometric agglomerated algebraic multigrid method with Gauss-Seidel
smoothing. The velocity and face fluxes are subsequently corrected using the updated
pressure field. Finally, the two turbulence transport equations are solved as implicit
advection-diffusion equations using an iterative linear solver with symmetric Gauss-Seidel smoothing, and the eddy viscosity is updated according to the
\(k\)-\(\omega\) SST model. This process is repeated until the residuals of the
velocity, pressure, and turbulence variables satisfy the prescribed convergence
criteria.

For the implementation, we use \texttt{OpenFOAM}
\cite{weller1998tensorial,OpenFOAM12}, an open-source finite volume CFD package widely
used in industrial applications. The
setup follows the \texttt{drivaerFastback} vehicle benchmark case.
The pressure correction equation is solved using \texttt{OpenFOAM}'s \texttt{GAMG} solver, whereas the momentum and turbulence equations are solved using \texttt{smoothSolver}.
For each sampled vehicle geometry, the surface mesh is obtained from~\texttt{DrivAerNet++}, and a geometry-specific volumetric mesh is generated using
the \texttt{OpenFOAM} meshing pipeline: a background mesh is first created, then locally refined,
and finally fitted to the vehicle surface using \texttt{snappyHexMesh}~\cite{OpenFOAMUserGuideSnappyHexMesh}.
The pressure coefficient fields used to train and evaluate the neural operator are
taken from the \texttt{DrivAerNet++} dataset~\cite{elrefaie2024drivaernet++}.
For the classical numerical solver study, we follow the simulation setup described in~\cite{elrefaie2024drivaernet++} but independently
compute geometry-specific reference solutions using \texttt{OpenFOAM}~12
on large volumetric meshes containing more than 20 million cells.
Within the class of
steady-RANS methods, this benchmark configuration is representative of current high-fidelity industrial practice in automotive aerodynamics.
Each reference RANS simulation is advanced for 7000 SIMPLE iterations.
Convergence is assessed using the residual histories and the stabilization
of the aerodynamic force coefficients over a moving window of 1000 iterations.

The cost of the implicit finite volume solver is the product of the number
of SIMPLE iterations $n_t$ and the cost per iteration.
Let \(n_c\) denote the number of control volumes and \(n_f\) the number of internal
faces in the finite volume mesh. For shape-regular polyhedral meshes dominated by hexahedral cells, one typically has \(n_f \approx 3n_c\).
Each SIMPLE iteration consists of assembling or updating the finite volume
operators and solving the resulting pressure, velocity, and turbulence model linear systems. Assuming each linear system takes about $m$ iterations of an
iterative solver, a rough cost model is 
\[
C_{\rm solver}
\approx
n_t \bigl(c_{\mathrm{assemble}} n_c
+ m c_{\rm solve} n_c
\bigr).
\]
The detailed values of \(c_{\mathrm{assemble}}\) and \(c_{\mathrm{solve}}\) are on the
order of a few thousand and a few hundred, respectively, and are reported in the supplement in
\ref{supp-ssec:fluid}.
The accuracy of the RANS finite volume solver depends not only on discretization error, but also on
mesh quality, geometric complexity, linear solver convergence, and fluctuations associated with turbulence modeling. In sufficiently smooth
regions and on nearly uniform meshes, second-order finite volume schemes typically yield
\[
\|U - U_h\|_{L^1} = \mathcal O(h^2),
\]
where $U$ and $U_h$ denote the exact and numerical solutions, respectively. 
On unstructured meshes for complex geometries, however, the effective convergence rate may deteriorate.
Assuming an empirical convergence rate
\[
\|U - U_h\|_{L^1} = \mathcal O(h^q),
\qquad q \leq 2,
\]
and sufficient trace regularity, we assume the relative $L^1$ error
in $C_p$ on the vehicle surface to inherit the same empirical rate.
Using \(n_c = \bigO(h^{-3})\) in three dimensions, the cost required to achieve an
error of order \(\varepsilon\) scales as
\begin{equation}
    C_{\rm solver}(\varepsilon)
=
\bigO\!\left(n_t m \varepsilon^{-3/q}\right)
\end{equation}
up to the dependence on the SIMPLE and inner linear solver iteration counts.

\subsubsection{Measured trade-off and acceleration mechanism}
We measure accuracy using the relative \(L^1\) error in the surface pressure coefficient, since integrated pressure loads depend directly on surface
integrals. For the classical baseline, we adapt the official \texttt{OpenFOAM}~12 \texttt{drivaerFastback} workflow~\cite{weller1998tensorial,OpenFOAM12} to six geometries sampled from \texttt{DrivAerNet++} and chosen to cover all vehicle categories.  For each geometry, we consider small, medium, and large volumetric meshes and vary the number of SIMPLE iterations. The corresponding 7000-iteration large-mesh solution serves
as the numerical reference. For the wall-clock measurements, the CPU allocation is scaled
with mesh size to maintain approximately \(10^5\) cells per core. The reported
times are therefore time-to-solution measurements under varying resource
allocations, rather than equal-core comparisons. 
\Cref{tab:openfoam_setup} reports the configurations and averaged mesh sizes, runtimes, and errors for the \texttt{OpenFOAM} predictions over six test geometries.  The canonical \texttt{drivaerFastback} configuration is reported separately in the supplement in \cref{supp-ssec:openfoam-fastback} as a representative individual case. 
The neural operator uses the unstructured-grid implementation. The input comprises the surface coordinates and postprocessed normals, giving
\(d_a=6\) channels on \(n_e\) points. Each layer incorporates the surface
normals through the geometric encoding proposed in~\cite{han2026geometric}, motivated
by geometry-dependent kernels such as double-layer
potentials~\cite{hess1967calculation}. 
We train on \(N=4000\) \texttt{DrivAerNet++} samples using a relative \(L^1\) loss. We fix \(k_{\max}=16\), \(d_g=64\), and \(L=4\), and vary the downsampled surface resolution over approximately \(10{,}000\), \(20{,}000\), and \(40{,}000\) points.

\begin{table}[htbp]
\centering
\begin{tabular}{l|cccc}
\Xhline{1.1pt}
Configuration
    & Large reference & Large & Medium & Small \\
\hline
Mean control volumes \(\overline{n_c}\)
    & \(22,401,336\) & \(22,401,336\) & \(2,962,931\)  & \(441,183\) \\
SIMPLE iterations $n_t$
    & 7000 & 2000 & 1000 & 1000 \\
CPU cores
    & 256 & 256 & 32 & 4 \\
\hline
Mean runtime (s)
    & \(6443.8\) & \(1860.1\) & \(921.8\) & \(547.0\) \\
Mean relative \(L^1\) error
    & Reference & \(12.06\%\) & \(18.20\%\) & \(25.89\%\) \\
\Xhline{1.1pt}
\end{tabular}
\caption{
\texttt{OpenFOAM} configurations and results for six selected vehicle geometries
drawn from \texttt{DrivAerNet++}. Mesh sizes, wall-clock runtimes, and
relative \(L^1\) errors in the surface pressure coefficient \(C_p\) are
averaged over test samples. For each geometry, the error is measured
relative to the corresponding geometry-specific reference solution
computed on the large mesh using \texttt{OpenFOAM}~12 and
7000 SIMPLE iterations.
}
\label{tab:openfoam_setup}
\end{table}

At approximately \(10{,}000\), \(20{,}000\), and \(40{,}000\) surface
points, the  average neural operator relative \(L^1\) errors are \(13.50\%\), \(12.66\%\), and \(12.03\%\), respectively. Increasing the surface resolution improves the representation of local geometric and flow features, but the gain is modest.  \Cref{fig:car-map} shows that the neural operator model reproduces the dominant pressure pattern for a representative case, qualitatively comparable to the \texttt{OpenFOAM} predictions on the medium mesh. 
For this dataset, neural operator accuracy appears to be limited in part by variability in the RANS-generated training data~\cite{wu2024transolver,alkin2025ab}.
In published benchmarks for \texttt{DrivAerNet++}, neural operator errors remain of order $10\%$~\cite{elrefaie2024drivaernet++,wu2024transolver,alkin2025ab}.

\begin{figure}[htbp]
     \centering
     \includegraphics[width=0.9\textwidth]{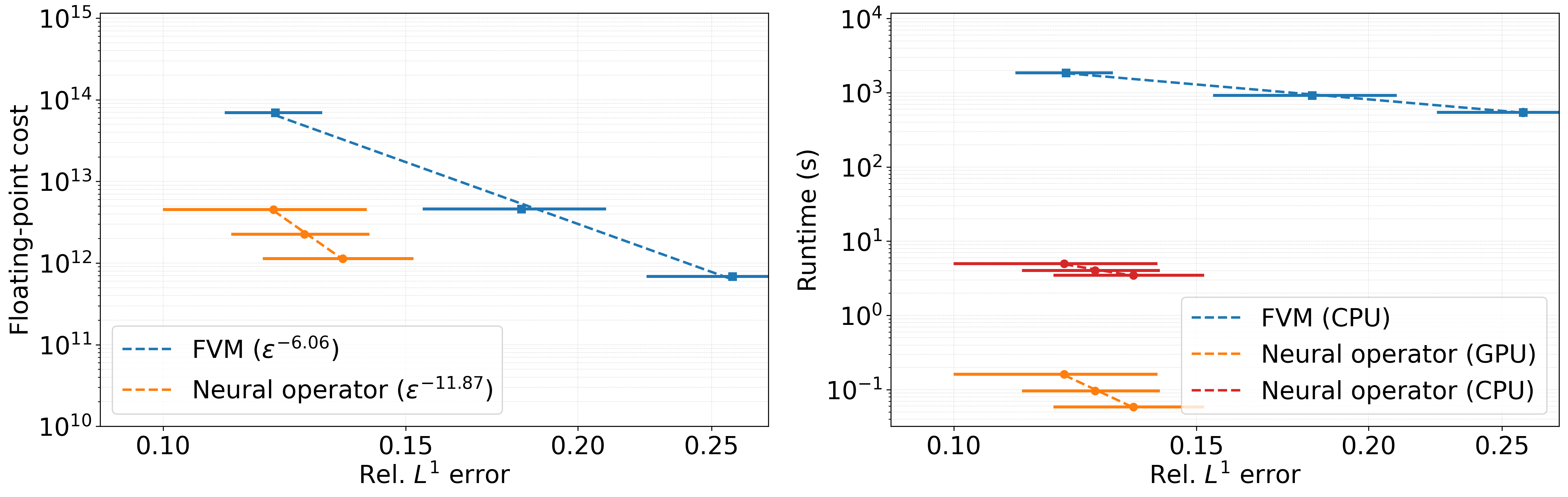}\\
     \caption{Cost-accuracy comparison for vehicle aerodynamics between the finite volume method (FVM) and the neural operator.
     Left: floating-point operation count versus relative $L^1$ error. 
     Right: wall-clock runtime in seconds versus relative $L^1$ error for CPU and GPU 
     implementations. Error bars denote one standard deviation over the test samples.
     }
     \label{fig:car-cost-accuracy}
\end{figure}

The cost-accuracy comparison is shown in~\cref{fig:car-cost-accuracy}. 
The result is consistent with two sources of online advantage. First, the
neural operator predicts the requested surface field on \(n_e\) points instead of computing a volumetric flow on \(n_c\) control volumes, with \(n_e\ll n_c\). Second, inference avoids the outer SIMPLE iteration and its associated linear solves. Its floating-point count is consequently below that of the large mesh FVM configuration, although it remains above that of the small mesh configuration. However, the steeper empirical slope of the neural operator curve suggests that this advantage may diminish in a higher-accuracy regime, even if the neural operator can attain that regime. 
GPU execution provides an additional wall-clock advantage.

\subsubsection{Limitations and implications}
\begin{figure}[htbp]
     \centering
     \includegraphics[width=0.9\textwidth]{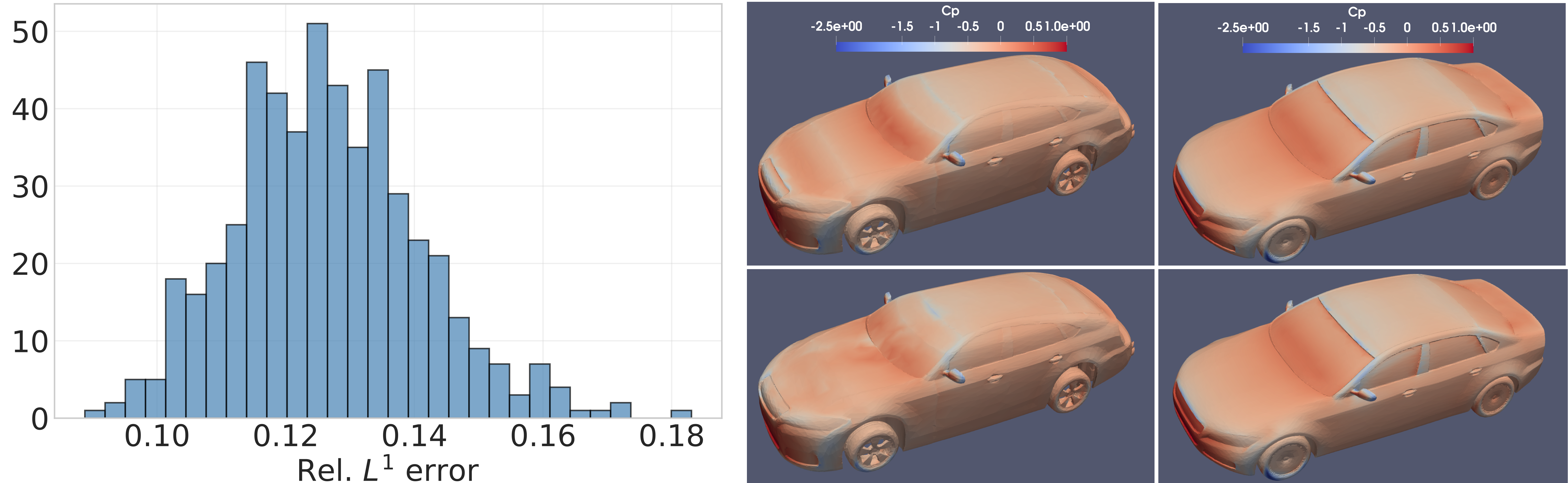}\\
     \caption{Error analysis for vehicle aerodynamics. 
     From left to right: the distribution of test errors, the test case with the largest error, and a representative median error case. For the two visualized cases, the reference solution is shown in the top row and the neural operator prediction in the bottom row.}
     \label{fig:car-error}
\end{figure}

The principal limitation is reliability across geometries.
Although the mean relative $L^1$ error is moderate, the distribution in \cref{fig:car-error} has a substantial upper tail. The largest error case corresponds to a relatively tall vehicle with strong pressure variations over the hood. Such prediction failures could compromise downstream engineering tasks,
including design optimization.
Moreover, because the surrogate predicts
only a surface quantity of interest, verification through a full PDE residual is not readily available. Reliable deployment will therefore
require uncertainty estimates and out-of-distribution indicators that identify questionable predictions and trigger verification, correction,
or replacement by a classical solver.

\section{Conclusion}
\label{sec:Conclusion}

We have studied when neural operator inference provides a favorable per-query cost-accuracy trade-off relative to strong classical solvers in the post-training, many-query limit. The experiments show that neural operators are competitive when they can attain the requested accuracy and either execute efficiently on the target hardware or bypass a substantial component of the classical computation. Across the cases considered, the observed benefits arise from hardware-efficient tensor operations, reduced sequential or iterative computation, and direct prediction of task-specific outputs. These mechanisms also explain why no method is uniformly superior: classical solvers reach smaller errors than those attained by the tested neural operators and remain essential for high-accuracy computation and verification.

The many-query limit is motivated by the prospect of scientific models that are reused across increasingly broad classes of problems. Although the neural operators studied here are problem-specific rather than
foundation models, the identified acceleration mechanisms indicate how broader reuse could increase their practical value. Realizing this potential will require reliable generalization, indicators of inaccurate or out-of-distribution predictions, and hybrid workflows in which classical solvers verify, correct, or replace neural predictions when higher accuracy or reliability is required.

\section*{Acknowledgments}
 DZH acknowledges the support of the high-performance computing platform of Peking University. After completing the initial draft, DZH used GPT-5.6 and DeepSeek-V3 to obtain suggestions for improving the language and presentation of the manuscript and the documentation of the accompanying software library. The authors reviewed all suggestions and incorporated them as appropriate; they take full responsibility for the final content. The authors also thank George Karniadakis for his helpful suggestions regarding this work.


\clearpage

\headers{Supplementary Materials: Cost-Accuracy Trade-offs}
        {D. Z. Huang, A. M. Stuart}

\begin{center}
{\bfseries
\MakeUppercase{Supplementary Materials}
\par}
\end{center}

\phantomsection
\addcontentsline{toc}{section}{Supplementary Materials}


\setcounter{section}{0}
\setcounter{subsection}{0}
\setcounter{subsubsection}{0}
\setcounter{equation}{0}
\setcounter{figure}{0}
\setcounter{table}{0}
\setcounter{theorem}{0}
\setcounter{algorithm}{0}

\renewcommand{\thesection}{SM\arabic{section}}
\renewcommand{\thefigure}{SM\arabic{figure}}
\renewcommand{\thetable}{SM\arabic{table}}

\section{Neural Operator}
\label{supp-sec:NO}
In this section, we give details of the neural operator architecture benchmarked in the present work. It is a modification of the Fourier neural operator~\cite{li2020fourier}, designed to better capture local effects~\cite{liu2024neural-gl,zeng2025point} and to accommodate problems with variable geometry~\cite{lingsch2023beyond,zeng2025point}. The neural operator was introduced and studied in \cite{zeng2025point,han2026geometric} under the names point cloud neural operator (PCNO) and multiscale point cloud neural operator (M-PCNO).  We describe the architecture and implementation in \cref{supp-ssec:nn_architecture,supp-ssec:nn_implementation}, and discuss its general cost-accuracy behavior in \cref{supp-ssec:cost-accuracy}.

\subsection{Neural operator architecture}
\label{supp-ssec:nn_architecture}
Given an input pair \((a,D)\), with \(a:D\to\mathbb R^{d_a}\), the neural operator defines a mapping between functions on the domain \(D\)
\begin{equation}
\label{supp-eq:target-operator}
    \omap:(a,D)\longmapsto u_{\rm q}.
\end{equation}

It first lifts the input \(a\) to a latent feature field \(g_0:D\to\mathbb R^{d_g}\), then applies \(L\) neural operator layers, and finally projects the resulting feature field onto the desired output space:
\begin{equation}
\label{supp-eq:NO_architecture}
\begin{split}
&g_0  = \cP\bigl(a\bigr), \qquad
g_{i} = \cL_i^D(g_{i-1}), \quad i=1,\dots,L, \qquad
u_{\rm q} =\cQ\bigl(g_L \bigr),\\
&\cG(a, D ;\theta) = \cQ\circ \cL^D_L \circ \cL^D_{L-1} \circ \cdots \circ \cL^D_1 \circ \cP (a).
\end{split}
\end{equation}
Here \(\cP:\mathbb R^{d_a}\to\mathbb R^{d_g}\) is a pointwise lifting map, \(\cQ:\mathbb
R^{d_g}\to\mathbb R^{d_u}\) is a pointwise projection map, and each \(\cL_i^D\)
is a neural operator layer acting on functions defined on \(D\) with $d_g$ channels. All hidden layers are taken to have the same channel dimension
\(d_g\). The lifting map increases the number of channels from \(d_a\) to \(d_g\), thereby embedding the input into a higher dimensional latent space. This higher dimensional representation can capture more complex features and make the underlying relationships more tractable~\cite{koopman1931hamiltonian,mezic2005spectral,huang2024operator}.

Suppressing the layer index, each neural operator layer $\cL^D$ maps the function space $\{g: D \rightarrow \R^{d_{g}}\}$ into itself. It consists of a linear operator followed by a pointwise nonlinear activation function, both defined at the operator level and
independent of the discretization.  The linear component is decomposed into a global operator
\(\mathcal K_{\rm long}\) and a local operator \(\mathcal K_{\rm short}\), so that 
\begin{equation}
\begin{split}
\label{supp-eq:mno}
    \cL^D g =  g + \sigma \Bigl(\mathcal{K}_{\rm long} g + \mathcal{K}_{\rm short} g \Bigr).
\end{split}
\end{equation}
Here, \(\mathcal K_{\rm long}\) captures long-range interactions, while
\(\mathcal K_{\rm short}\) captures short-range effects. The activation function
\(\sigma\), such as \texttt{GELU}~\cite{hendrycks2016gaussian}, acts pointwise at each
spatial location and on each channel. 
The residual connection is inspired by residual networks~\cite{he2016deep} and the ODE interpretation of deep neural networks~\cite{WeinanE20171}, and it helps stabilize training, especially for deeper architectures.

For the global component of the linear operator, we use a Fourier neural operator (FNO) representation~\cite{li2020fourier,kovachki2023neural,nelsen2021random,nelsen2024operator,kossaifi2024library}.
The corresponding global kernel is represented by a truncated Fourier expansion:
\begin{equation}
\label{supp-eq:FNO-kernel}
(\mathcal K_{\rm long} g)(x)
=
\sum_{\|k\|_\infty\le k_{\max}}
\int_{D}\kappa_{k}(x,y)W_k^v\,g(y)\,dy 
=
\sum_{\|k\|_\infty\le k_{\max}}
\phi_k(x)\,
W_k^v
\int_{D}\overline{\phi_k(y)}\,g(y)\,dy,
\end{equation}
where \(k\in\mathbb Z^d\) denotes a Fourier mode, $\kappa_k(x,y)=
e^{2\pi i k\cdot \frac{x-y}{l}}$, $\phi_k(x)=e^{2\pi i k\cdot x/l}$,
\(l\) is a characteristic length scale of the domain, and
\(W_k^v\in\mathbb C^{d_g\times d_g}\) is the learned complex matrix associated with mode
\(k\). 
Only modes satisfying \(\|k\|_\infty\le k_{\max}\) are retained, while all remaining modes are
set to zero. This truncation reduces the number of learnable parameters and yields a
low-rank spectral representation of the global operator.
When the computational domain $D$ is a fixed periodic box discretized on a structured grid,
\eqref{supp-eq:FNO-kernel} can be evaluated efficiently using the fast Fourier transform (FFT):
\begin{align}
\label{supp-eq:FNO}
    \mathcal{K}_{\rm long} g = \cF^{-1} \Bigl(R \cdot \bigl(\cF g\bigr)\Bigr), 
\end{align}
where \(\mathcal F\) and \(\mathcal F^{-1}\) denote the channelwise Fourier transform
and inverse Fourier transform, respectively. The multiplier \(R\) is defined by
\[
R(k)=W_k^v \quad \text{for } \|k\|_\infty\le k_{\max},
\qquad
R(k)=0 \quad \text{otherwise}.
\]
Thus, the FNO layer applies a learned linear map to each retained Fourier mode and then
transforms the result back to physical space.
Alternative kernel-based neural operators use kernels that are localized in physical space~\cite{li2020neural}, for example satisfying
\(\kappa_k(x,y)=0\) whenever \(\lVert x-y\rVert_2>r\) for a prescribed interaction radius \(r>0\), or kernels that are made input-dependent through attention mechanisms~\cite{cao2021choose,calvello2025continuum}. The Fourier basis \(\{\phi_k\}\) may also be replaced by multiwavelet bases~\cite{gupta2021multiwavelet} or Laplace-Beltrami eigenbases~\cite{bonev2023spherical,chen2024learning}.

For the local component of the linear operator, we use a pointwise affine map together with a gradient correction~\cite{zeng2025point,liu2024neural-gl}:
\begin{equation}
\begin{split}
    \label{supp-eq:local-operator}
    &(\mathcal{K}_{\rm short} g)(x)  =  
    W g(x) + b + W^{g} \textrm{vec}\bigl(\nabla_D g(x)\bigr),
\end{split}
\end{equation}
where \(W\in\mathbb R^{d_g\times d_g}\), \(b\in\mathbb R^{d_g}\), and
\(W^g\in\mathbb R^{d_g\times d_g d}\) are learned parameters. 
The pointwise affine term is important because integral operators with continuous kernels are compact on many
infinite-dimensional function spaces, whereas the identity operator is not \cite[Section 2]{kress1989linear}. Consequently, a purely integral architecture may not represent identity-type components efficiently. The gradient term is
motivated by the structure of PDEs, where differential operators and their compositions
appear naturally. Including \(\nabla_D g\) therefore helps the model capture local
high-frequency features and sharp spatial variations.
Near discontinuities, discrete gradients may become large and lead to numerical
instability. To regularize this effect, we apply the \texttt{SoftSign} activation componentwise, followed by a learned scaling matrix $W^{g'}\in\R^{d_g\times d_g}$:
\begin{equation}
\label{supp-eq:smoothed-gradient}
    \mathcal{K}_{\rm short} g(x)  =  
    W g(x) + b + W^{g'}\texttt{SoftSign}\Bigl(W^{g}\textrm{vec}\bigl(\nabla_D g(x)\bigr) \Bigr) \ \textrm{ where }\  \texttt{SoftSign}(z) = \frac{z}{1 + |z|}.
\end{equation}
The smoothed gradient correction in~\eqref{supp-eq:smoothed-gradient} preserves gradient information while also serving as an indicator of sharp or discontinuous local features in $g$.

Combining \eqref{supp-eq:FNO-kernel}, or equivalently \eqref{supp-eq:FNO}, with \eqref{supp-eq:smoothed-gradient} defines the linear component of the neural operator layer. Both components depend implicitly on the computational domain \(D\). The learnable parameters in each layer are collected in $\{W_k^v, W, b, W^g, W^{g'}\}$. The resulting operator is designed at the continuous level and subsequently discretized according to the underlying geometry.   We discuss the implementation details in the next section.

\subsection{Implementation details}
\label{supp-ssec:nn_implementation}
The neural operator requires the evaluation of both a global integral operator and a local gradient operator. The corresponding implementation depends on the geometry of the computational domain $D$.

\paragraph{Structured grids}
When the computational domain $D$ is a box discretized on a structured grid, the global integral operator is evaluated using the fast Fourier transform, as in \eqref{supp-eq:FNO}, while the local gradient is approximated using centered finite differences. Specifically, for each spatial direction \(i=1,\dots,d\),
\begin{equation}
\label{supp-eq:central-diff}
    \frac{\partial g}{\partial x_i}(x) = \frac{g(x + h_i  e_i) - g(x  - h_i e_i)}{2h_i} ,
\end{equation}
where $h_i$ is the mesh spacing in the \(i\)-th direction and \(e_i\) is the
\(i\)-th Euclidean basis vector. At grid points adjacent to the boundary, the centered difference is evaluated using periodic wraparound indexing.

\paragraph{Unstructured grids and variable geometries}
When the computational domain is general or variable, we discretize it using a point cloud $\{x^{(i)}\}_{i=1}^{n_e} \subset D$, together with mesh connectivity information and quadrature weights
\(\{\dd\Omega_i\}_{i=1}^{n_e}\). In this setting, the global integral operator is evaluated by applying a separable quadrature rule to \cref{supp-eq:FNO-kernel}:
\begin{equation}
\label{supp-eq:unstructured-integral}
\begin{split}
    &(\mathcal K_{\rm long} g)(x) \approx 
    \sum_{\|k\|_\infty\le k_{\max}}
\phi_k(x)\,
W_k^v
\sum_{i=1}^{n_e}\overline{\phi_k(x^{(i)})}\,g(x^{(i)})\,\dd\Omega_i.
\end{split}
\end{equation}
Thus, for each retained mode \(k\), one first computes the weighted modal coefficient, then applies the learned matrix \(W_k^v\), and finally reconstructs the output
at the target point \(x\).
For the local gradient, we follow the point-cloud least-squares construction~\cite{zeng2025point}. For each point
\(x\), let 
\(\{x^{(j)}\}_{j=1}^{\nu(x)}\) denote its neighboring points, where $\nu(x)$ is the number of neighbors. The discrete gradient
\(\nabla_D g(x)\in\mathbb R^{d_g\times d}\) is defined as the solution of the
local least-squares problem
\begin{equation}
\label{supp-eq:pcno-gradient-ls}
\nabla_D g(x) 
\approx
\arg\min_{G \in\mathbb R^{d_g \times d}}
\sum_{j=1}^{\nu(x)}
\left|
g(x^{(j)})-g(x)- G\cdot(x^{(j)}-x)
\right|^2.
\end{equation}
We introduce the geometry matrix and the feature-difference matrix
\begin{align}
&A(x) = \begin{bmatrix}
x^{(1)} - x\quad \cdots \quad x^{(\nu(x))} - x
\end{bmatrix} \in \R^{d \times \nu(x)},
\\ 
&B_g(x) = \begin{bmatrix}
g(x^{(1)}) - g(x) \quad 
\cdots \quad 
g(x^{(\nu(x))}) - g(x)
\end{bmatrix} \in \R^{d_g \times \nu(x)}.
\end{align}
Then \eqref{supp-eq:pcno-gradient-ls} can be written equivalently as
\begin{equation}
\nabla_D g(x) 
\approx
\arg\min_{G \in\mathbb R^{d_g \times d}} \lVert G A(x)  -  B_g(x) \rVert_F^2
\end{equation}
and the least-squares solution is
\begin{equation}
\label{supp-eq:pcno-gradient-solution}
\nabla_D g(x) = B_g(x)\,A(x)^\dagger,
\end{equation}
where $A(x)^{\dagger} \in \R^{\nu(x) \times d}$ is the pseudoinverse of \(A(x)\),
precomputed from the local geometry. If the $j$-th row of $A(x)^{\dagger}$ is denoted by $A_j(x)^{\dagger} \in \R^{d}$, then \cref{supp-eq:pcno-gradient-solution} can be evaluated in a single message-passing step:
\begin{equation}
\label{supp-eq:gradient_operator}
    \nabla_D g(x) = \sum_{j=1}^{\nu(x)} \bigl(g(x^{(j)}) - g(x)\bigr) A_j(x)^{\dagger}.
\end{equation}
When \(D\) is a surface embedded in \(\mathbb R^d\), as in the vehicle aerodynamics case, this procedure computes a least-squares gradient in ambient coordinates, which gives a tangential gradient approximation on the surface.

\begin{remark}[Boundary conditions]
\label{supp-rem:bc}
In the architectures considered here, boundary conditions are generally not imposed strongly by the neural operator layers themselves, but are instead learned from the training data.
Periodic boundary conditions are an exception, since they are naturally compatible with
the Fourier representation in \eqref{supp-eq:FNO} and with periodic wraparound indexing.
For nonperiodic box domains, we follow the padding strategy used in the FNO~\cite{li2020fourier}: the domain is enlarged by \(10\%\) in each coordinate
direction and the data are zero-padded to provide a simple periodic extension in the spirit
of Fourier continuation~\cite{bruno2007accurate}. For unstructured grids and variable
geometries, the characteristic length $l$ in the Fourier features is chosen to be approximately twice the size of the bounding box, following~\cite{han2026geometric}. 
Boundary conditions may also be imposed strongly, for example by using bases adapted to prescribed boundary conditions
\cite{hesthaven2018non,bhattacharya2021model,liu2023spfno,chen2024learning,stuart2025enforcing}. In the present  work, we do not consider problems with variable boundary conditions. When such
conditions vary across samples, they can be incorporated into the input data \(a\), for
example through additional boundary channels; see \cite{zhou2025boundary}.
\end{remark}

\subsection{Cost-accuracy behavior}
\label{supp-ssec:cost-accuracy}
Finally, we discuss the cost-accuracy behavior of the neural operator described above. We focus on online inference cost. Floating-point operation counts are estimated under the standard model in which each multiplication or addition is counted as one flop, while memory traffic, communication, and parallelization overhead are ignored.

Let \(n_e\) denote the total number of grid points or point-cloud degrees of freedom, and let
\[
K := (2k_{\max}+1)^d
\]
denote the number of retained modes in the truncated spectral representation. Recall that \(d_g\) and \(L\) denote the hidden channel dimension and the number of
neural operator layers, respectively. Together, \(n_e\), \(K\), \(d_g\), and \(L\) determine the leading-order inference cost.
The precise cost depends on the geometry of the computational domain.
On a structured periodic grid, the global operator is evaluated by FFT, whereas on an unstructured point cloud it is evaluated by the separable quadrature formula \eqref{supp-eq:unstructured-integral}.
The local gradient correction is computed by centered finite differences on structured grids and by least-squares message passing on unstructured grids. Since the local pseudoinverses $A(x)^{\dagger}$ are precomputed, the two gradient implementations have comparable leading-order complexity.

\paragraph{Lifting and projection}

Assume that the lifting map $\cP$ is a
pointwise affine map. Then its cost is
\[
C_{\rm lift} = 2n_e d_a d_g.
\]
Assume further that the projection map $\cQ$ is a two-layer pointwise multilayer perceptron $\mathbb R^{d_g}\to\mathbb R^{d_g}\to\mathbb R^{d_u}$ with one pointwise activation in the hidden layer. Then its cost is
\[
C_{\rm proj} = 2n_e d_g^2 + 2n_e d_g d_u + c_{\sigma} n_e d_g.
\]

\paragraph{Structured-grid cost}
For a structured grid, the global integral operator is evaluated using the FFT.
If one
complex FFT of size \(n_e\) costs \(5n_e\log_2 n_e\) flops, then one forward and one
inverse transform on all \(d_g\) channels cost  $10 d_g n_e \log_2 n_e$.
For the retained modes, each complex multiplication
is counted as 6 flops, hence each complex matrix-vector product costs
$d_g(6d_g+2(d_g-1)) = 8d_g^2-2d_g$
flops. Hence the total cost of the structured global operator is
\begin{equation}
10 d_g n_e \log_2 n_e
+
K(8d_g^2-2d_g).
\end{equation}
For the local operator \eqref{supp-eq:local-operator}, the pointwise affine term costs $2d_g^2 n_e$.
The centered finite difference gradient \eqref{supp-eq:central-diff} requires 
one subtraction and one multiplication per point, per channel, and per spatial direction, and therefore costs $2 d d_g n_e$.
Applying the gradient map \(W^g  \in \R^{d_g \times (d_g d)}\) costs $(2d\,d_g^2-d_g)n_e$. 
Applying the componentwise \(\mathrm{SoftSign}\) nonlinearity costs
$3d_g n_e$.
Applying the scaling matrix $W^{g'}\in \R^{d_g \times d_g}$ costs $(2\,d_g^2-d_g)n_e$.
Adding the pointwise branch and the gradient-correction branch costs \(d_g n_e\), and
combining the local contribution with the global contribution costs another \(d_g n_e\).
Finally, evaluating the outer activation and applying the residual connection cost $c_\sigma d_g n_e + d_g n_e.$ 

Therefore, one structured-grid neural operator layer
costs
\begin{equation}
\label{supp-eq:cost-layer-structured}
\begin{aligned}
C_{\rm layer}
&=
10 d_g n_e \log_2 n_e
+
K(8d_g^2-2d_g)
+
((2d +4) d_g^2 + (2d + 4 + c_{\sigma}) d_g)n_e. 
\end{aligned}
\end{equation}
For an \(L\)-layer structured-grid neural operator, the inference cost satisfies
\begin{equation}
\begin{split}
    C^{\rm infer}_{\rm NO} &= C_{\rm lift} + C_{\rm proj} + L C_{\rm layer}.
\end{split}
\end{equation}
When \(K\), \(n_e\), and \(d_g\) are the dominant parameters, this reduces to
\begin{equation}
\label{supp-eq:cost-structured}
\begin{split}
    C^{\rm infer}_{\rm NO} &= \bigO\bigl(10 L d_g n_e \log_2 n_e + 8 L K d_g^2 + (2 + (2d+4)L) d_g^2 n_e\bigr).
\end{split}
\end{equation}

\paragraph{Unstructured-grid cost}
For a general point cloud, the global integral operator is evaluated using
\eqref{supp-eq:unstructured-integral}. To reduce memory usage, especially in the case of
variable computational domains, we evaluate the basis functions $\phi_k(x)$ online rather than storing them for all modes and all points. For each retained mode \(k\),
evaluating the basis on all points costs $4 d n_e$ flops, forming the weighted modal coefficient costs \(4d_g n_e\) flops, applying the learned matrix \(W_k^v\) costs \(8d_g^2-2d_g\) flops, and reconstructing the output at all points costs \(8d_g n_e\) flops. 
Hence the total cost of the unstructured global operator is
\begin{equation}
\label{supp-eq:cost-unstructured-integral}
K (12d_g + 4d) n_e + K(8d_g^2 - 2d_g).
\end{equation}
For the local operator, the only difference from the structured-grid case lies in the
gradient evaluation.  In the
message-passing gradient formula \eqref{supp-eq:gradient_operator}, each directed neighbor
contribution requires \(d_g\) flops to form the feature difference and \(2d\,d_g\) flops
to apply the local aggregation. Thus, the gradient cost is
\begin{equation}
2(2d+1)d_g\,n_f.
\end{equation}
The remaining pointwise operations are the same as in the structured-grid case. The
pointwise affine term costs \(2d_g^2 n_e\); the gradient map \(W^g\) costs
\((2d\,d_g^2-d_g)n_e\); the componentwise \(\mathrm{SoftSign}\) costs
\(3 d_g n_e\); the scaling matrix \(W^{g'}\) costs \((2d_g^2-d_g)n_e\); the
branch additions cost \(2d_g n_e\); and the outer activation and residual connection cost
\(c_\sigma d_g n_e + d_g n_e\).

Therefore, one unstructured-grid neural operator layer costs
\footnote{When the surface normals are incorporated in each layer following the geometric encoding strategy of \cite{han2026geometric}, the global operator costs
\begin{equation*}
K (12d_g + 4d) n_e + K(8d_g^2 - 2d_g) + 2d_g(d + 2d_g (d+1) - 1) n_e.
\end{equation*}
The local operator and the remaining terms cost 
\begin{equation*}
2(2d+1)d_g\,n_f  
+
\bigl((2d+4)d_g^2 + (4+c_\sigma)d_g\bigr)n_e
+  d_g (4d_g + 2(d^2+d) + 2) n_e
\end{equation*}
The cost for computing the gradient of the surface normal is $2(2d+1) d n_f$. Hence, for an \(L\) layer unstructured-grid neural operator, the total inference cost up to low-order terms is
\begin{equation}
\begin{split}
    C^{\rm infer}_{\rm NO} &= C_{\rm lift} + C_{\rm proj} + L C_{\rm layer} \\
                 &= \bigO\bigl(12 K L d_g n_e + 8 L K d_g^2 + (2 + (6d+12)L) d_g^2 n_e\bigr).
\end{split}
\end{equation}
}
\begin{equation}
\begin{aligned}
C_{\rm layer} 
&=
K (12d_g + 4d) n_e + K(8d_g^2 - 2d_g)
+
2(2d+1)d_g\,n_f  
+
\bigl((2d+4)d_g^2 + (4+c_\sigma)d_g\bigr)n_e .
\end{aligned}
\end{equation}
In practice, the neighborhood size is typically uniformly bounded and is comparable to
that of a structured grid, so that  \(n_f \approx d n_e\).
Hence, for an \(L\) layer unstructured-grid neural operator, the total inference cost up to low-order terms is
\begin{equation}
\label{supp-eq:cost-unstructured}
\begin{split}
    C^{\rm infer}_{\rm NO} &= C_{\rm lift} + C_{\rm proj} + L C_{\rm layer} \\
                 &= \bigO\bigl(12 K L d_g n_e + 8 L K d_g^2 + (2 + (2d+4)L) d_g^2 n_e\bigr).
\end{split}
\end{equation}

For fixed \(K\), \(d_g\), and \(L\), the structured-grid cost \eqref{supp-eq:cost-structured} is quasilinear in \(n_e\), whereas the unstructured-grid cost \eqref{supp-eq:cost-unstructured}  is linear in \(n_e\). However, the associated prefactors may
nevertheless be substantial because they depend strongly on the hidden channel dimension, the number of retained modes, and the network depth.

\paragraph{Predictive accuracy}
The predictive accuracy of neural operators remains less well understood than that of classical numerical solvers. A useful heuristic model~\cite{lu2021learning,de2022cost} for the in-distribution test error is  
\begin{equation}
\label{supp-eq:empirical_predictive_error}
\E_{(a,D)\sim\mu}\Bigl( \frac{\|\omap(a,D)-\cG(a,D;\theta_N)\|_{\cU}}{\|\omap(a,D) \|_{\cU}}\Bigr) 
\approx
c_{\rm disc} h_{\rm eff}^q + c_{\rm data} N^{-\beta} + \varepsilon_{\rm opt},
\end{equation}
where \(\theta_N\) denotes parameters trained on \(N\) samples.
The terms in this decomposition correspond to distinct sources of error and have been analyzed only in simplified settings.
The first term represents discretization or representation error. Here \(h_{\rm eff}\) denotes an effective resolution scale. When \(k_{\max}\) is sufficiently large, the Fourier representation can resolve the underlying spatial discretization, and \(h_{\rm eff}\) may be identified with the physical mesh size~\cite[Theorem 24]{kovachki2021universal}\cite{kim2024bounding}. When \(k_{\max}\) is more restrictive, the dominant error may instead arise from Fourier truncation, in which case 
$h_{\rm eff} \approx k_{\max}^{-1}$ may be regarded as the effective resolution scale~\cite[Theorem 2.1]{han2026geometric}. Related representation errors also arise from encoder-decoder architectures, such as PCA-based neural operators \cite[Theorem~3.4]{bhattacharya2021model}.
The exponent \(q>0\) should therefore be interpreted as an effective approximation exponent, reflecting both the regularity of the target operator and the approximation properties of the neural operator architecture. 
The second term represents the statistical error arising from the approximation of the population risk
\begin{align}
\label{supp-eq:risk}
                      \cJ_{\infty}(\theta):=\E_{(a,D)\sim\mu} \Bigl( \frac{\|\omap(a,D)-\cG(a,D;\theta)\|_{\cU}}{\|\omap(a,D)\|_{\cU}} \Bigr)
\end{align}
by the empirical risk~\eqref{supp-eq:emp_risk} 
\begin{align}
\label{supp-eq:emp_risk}
                    \cJ_{N}(\theta):=\E_{(a,D)\sim\mu^N}\Bigl( \frac{\|\omap(a,D)-\cG(a,D;\theta)\|_{\cU}}{\|\omap(a,D) \|_{\cU}}\Bigr)  = \frac{1}{N}\sum_{i=1}^N\frac{\|\omap(a_i,D_i)-\cG(a_i,D_i;\theta)\|_{\cU}}{\|\omap(a_i,D_i)\|_{\cU}}
\end{align}
based on \(N\) training samples with $\mu^N = \frac{1}{N}\sum_{j=1}^N \delta_{(a_j, D_j)}$. 
The exponent \(\beta>0\) generally depends on the regularity of the target solution operator and on the input distribution; see, for example, \cite[Theorem 1.3]{de2023convergence}\cite[Corollary 2.8]{reinhardt2026statistical}\cite[Theorem 4.9]{huang2024operator}\cite{liu2024neural}\cite[Theorem 4]{liu2024deep}. 
The final term, \(\varepsilon_{\rm opt}\), accounts for residual error because
training only approximately minimizes the nonconvex empirical objective. This
residual may arise from finite optimization time, stochastic gradient noise, or
convergence to a suboptimal stationary point
\cite{bottou2007tradeoffs,cisneros2025optimization}.

Accordingly, in the present work, we focus primarily on empirical predictive accuracy rather than on a formal a priori error analysis.

\section{Classical numerical solvers}
\label{supp-ssec:classical_solvers}
In this section, we derive detailed floating-point operation counts for the classical numerical solvers used in the three benchmarks; comparable derivations are not readily available in the literature. \Cref{supp-ssec:darcy} considers finite elements with geometric multigrid for two-dimensional Darcy flow, \cref{supp-ssec:ns} considers a Fourier pseudospectral method with fourth-order Runge-Kutta time integration for two-dimensional incompressible Navier-Stokes flow, and \cref{supp-ssec:fluid} considers a finite-volume RANS solver for vehicle aerodynamics. Throughout, each multiplication or addition is counted as one floating-point operation, while memory traffic, communication, and parallelization overhead are excluded.

\subsection{Darcy flow}
\label{supp-ssec:darcy}
Assume that the multigrid hierarchy consists of \(L+1\) levels, indexed by
\(l=0,1,\dots,L\), where level \(0\) is the finest grid and level \(L\) is the
coarsest grid. We assume that each coarsening step reduces the number of elements in
each spatial direction by a factor of \(2\), and that the linear system on the coarsest
grid is solved directly.

At level \(l\), the domain is partitioned into an
\(n_{l}\times n_{l}\) quadrilateral mesh. For conforming \(Q_1\) finite
elements with homogeneous Dirichlet boundary conditions, the number of active degrees of
freedom is 
$$n_{e,l} = (n_{l}-1)^2.$$
Since each interior degree of freedom couples to the degrees of freedom in the
surrounding \(3\times 3\) nodal stencil, the number of nonzero entries in the stiffness
matrix \(A_l\) is
\[
\operatorname{nnz}(A_l)=(3n_{l}-5)^2.
\]
Let \(P_l\) denote the prolongation matrix from level \(l+1\) to level \(l\),
that is, from the coarse grid to the next finer grid. For bilinear interpolation on a
uniform quadrilateral grid, 
\[
\operatorname{nnz}(P_l)=9 n_{e,l+1}.
\]
The restriction matrix from  level \(l\) to level \(l+1\) is taken to be \(P_l^T\).

The stiffness matrices on all levels, the prolongation operators between adjacent
levels, and the finest-grid load vector are assembled in the assembly phase. The assembly cost therefore satisfies
\begin{equation}
  C_{\rm assemble} \approx \sum_{l = 0}^{L} 156 n_{l} ^2 + 60 n_{0}^2 + \sum_{l = 1}^{L} 9 n_{e,l}  \approx 271 n_{e}.
\end{equation}
Here \(156\) is the per-element cost for assembling the stiffness matrix, \(60\) is the
per-element cost for assembling the finest-grid load vector, and the third term bounds
the cost of constructing the prolongation operators. 
The constants \(156\) and \(60\) are obtained from a direct quadrature-level flop count for \(Q_1\) quadrilateral elements with a \(2\times2\) tensor-product Gaussian rule.
At each quadrature point, evaluating the nodal \(Q_1\) coefficient \(a_h\) requires
\(4\) multiplications and \(3\) additions, for a total of \(7\) flops. Updating the
\(4\times4\) local stiffness matrix requires one multiplication and one addition per
entry, namely \(32\) flops. Thus the stiffness assembly costs \(7+32=39\) flops per
quadrature point, and hence $4\times 39 = 156$
flops per element.
Similarly, evaluating the nodal \(Q_1\) source term \(f_h\) costs \(7\) flops per
quadrature point, while updating the four entries of the local load vector costs
\(8\) flops. Thus the load-vector assembly costs \(7+8=15\) flops per quadrature point,
and hence
$4\times 15 = 60$ flops per element. Only the finest-grid load vector is assembled from the source term;
on coarser levels, the right-hand sides are obtained by restricting residuals during the
multigrid iteration.

We next estimate the floating-point cost of one multigrid \(V\)-cycle. At each non-coarsest level
\(l=0,1,\dots,L-1\), the algorithm consists of \(\nu_1\) pre-smoothing steps,
one residual computation, one restriction to the next coarser level, one prolongation of
the coarse-grid correction, and \(\nu_2\) post-smoothing steps. Let
$\nu := \nu_1+\nu_2$
denote the total number of smoothing steps per level. 
For Richardson-Jacobi smoothing,
$$u_{h,l} \leftarrow  u_{h,l} + \omega D_{l}^{-1} (b_l - A_l u_{h,l}),$$
where $b_l$ is the current right-hand side (equal to the restricted residual), $D_l$ is the diagonal of $A_l$, and $\omega$ is the relaxation parameter.
One smoothing step costs $2\,\operatorname{nnz}(A_l) + 4 n_{e,l}$,
where \(2\,\operatorname{nnz}(A_l)\) is the cost of the sparse matrix-vector product,
and \(4n_{e,l}\) accounts for residual formation, diagonal scaling, and the vector update.
The residual is then computed as 
$$r_{l} = b_l - A_{l} u_{h,l},$$ 
which costs $2\,\operatorname{nnz}(A_l)+n_{e,l}$.
The restriction step 
$$b_{l+1} = P_l^T r_l$$ 
projects the residual onto the right-hand side at the next coarser level and costs $2\,\operatorname{nnz}(P_l)$ flops. 
The prolongation of the coarse-grid correction $e_{l+1}$ back to level $l$, followed by the correction step
$$u_{h,l} \leftarrow u_{h,l} + P_l e_{l+1}$$
costs $2\,\operatorname{nnz}(P_l)+ n_{e,l}$. 
On the coarsest level, the linear system is solved directly by LU factorization.
The corresponding cost is approximately \(n_{e,L}^3\), which is \(\bigO(1)\) when the coarsest grid is fixed.
Therefore, the cost of one \(V\)-cycle is
\begin{equation}
\begin{split}
C_{\mathrm{solve}}
&\approx
\sum_{l=0}^{L-1}\Bigl(\nu\bigl(2\,\operatorname{nnz}(A_l)+ 4 n_{e,l}\bigr)
+\bigl(2\,\operatorname{nnz}(A_l)+n_{e,l}\bigr)
+\bigl(4\,\operatorname{nnz}(P_{l})+n_{e,l}\bigr)\Bigr) + n_{e,L}^3 \\
&\approx (\frac{88}{3}\nu + \frac{116}{3})  n_{e} .
\end{split}
\end{equation}

Combining the assembly cost with the solve cost of \(m\) multigrid \(V\)-cycles yields
\begin{equation}
C_{\rm solver}
=
C_{\mathrm{assemble}}+mC_{\mathrm{solve}}
\approx 271 n_{e} + m (\frac{88}{3}\nu + \frac{116}{3})  n_{e} = \left(
c_{\rm assemble}
+
c_{\rm solve}\, m
\right) n_e.
\end{equation}
Therefore, $c_{\rm assemble} = 271$ and $c_{\rm solve} = \frac{88}{3}\nu + \frac{116}{3}$. For the common choice $\nu_1 = \nu_2 = 3$,  both coefficients are on the order of a few hundred.

\subsection{Incompressible Navier-Stokes flow}
\label{supp-ssec:ns}
Assume that a two-dimensional complex FFT on an \(n_e = n\times n\) grid 
costs $5 n_e \log_2 n_e$.
Given the Fourier coefficients \(\widehat{\omega}\), the right-hand side is evaluated in spectral form as
\[
\widehat{\mathcal N}(\widehat{\omega})_{k}
=
-\widehat{(v\cdot\nabla \omega)}_{k}
-\nu |2\pi k|^2 \widehat{\omega}_{k}
+\widehat f_{k},
\]
where \(k\in\mathbb Z^2\) is the integer Fourier index and \(2\pi k\)
is the corresponding physical wavenumber on \(D=[0,1]^2\).
The nonlinear term \(v\cdot\nabla\omega\) is evaluated pseudo-spectrally: \(v\) and
\(\nabla\omega\) are transformed to physical space, multiplied pointwise, and then
transformed back to Fourier space.
For  $|k|\neq 0$, the velocity and vorticity derivatives are computed in
Fourier space as
\[
\widehat v_1 = \frac{i 2\pi k_2}{|2\pi k|^2}\widehat\omega,
\qquad
\widehat v_2 = -\frac{i 2\pi k_1}{|2\pi k|^2}\widehat\omega,
\qquad
\widehat{\partial_x\omega}=i 2\pi k_1\widehat\omega,
\qquad
\widehat{\partial_y\omega}=i 2\pi k_2\widehat\omega .
\]
These four Fourier multiplier operations cost approximately \(24 n_e\) flops if each
complex multiplication is counted as \(6\) flops. We then apply four inverse FFTs to
obtain \(v_1,v_2,\partial_x\omega,\partial_y\omega\) in physical space. The nonlinear
term
\[
v\cdot\nabla\omega
=
v_1\partial_x\omega+v_2\partial_y\omega
\]
is formed pointwise at a cost of \(3 n_e\) flops. A forward FFT is then used to transform the nonlinear term
back to Fourier space. 
Finally, the viscous term \(-\nu |2\pi k|^2\widehat\omega\), the
forcing term, and the nonlinear term are combined in Fourier space, at a cost of
approximately \(10n_e\) flops. 
Dealiasing using the $2/3$ rule~\cite{orszag1972numerical} is applied at each right-hand-side
evaluation; for a simple spectral mask, this costs approximately \(2n_e\) flops.
Therefore, one right-hand-side evaluation costs
\begin{align}
C_{\mathrm{rhs}}
&\approx
25 n_e \log_2 n_e
+(24+3+10+2)n_e = 25 n_e \log_2 n_e
+ 39 n_e. 
\end{align}
One RK4 time step requires four right-hand-side evaluations,
\[
\begin{aligned}
\widehat{k}_1 &= \widehat{\mathcal N}(\widehat{\omega}^n),\quad 
\widehat{k}_2 &= \widehat{\mathcal N}(\widehat{\omega}^n+\tfrac{\Delta t}{2}\widehat{k}_1),\quad 
\widehat{k}_3 &= \widehat{\mathcal N}(\widehat{\omega}^n+\tfrac{\Delta t}{2}\widehat{k}_2),\quad 
\widehat{k}_4 &= \widehat{\mathcal N}(\widehat{\omega}^n+\Delta t\,\widehat{k}_3),
\end{aligned}
\]
followed by the spectral-space update
\[
\widehat{\omega}^{n+1}
=
\widehat{\omega}^n+
\frac{\Delta t}{6}(\widehat{k}_1+2\widehat{k}_2+2\widehat{k}_3+\widehat{k}_4).
\]
The RK4 stage-state formation and final update require 
approximately \(28n_e\) additional
flops. 
If $n_t$ time steps are used to reach the final time \(T\), the total
floating-point cost is
\begin{equation}
\begin{split}
    C_{\rm solver}
&\approx n_t(4C_{\mathrm{rhs}} + 28n_e) = 
n_t\left(100 n_e\log_2 n_e+184 n_e\right)  \\
&= n_t\left(c_{\rm assemble}\, n_e\log_2 n_e + c'_{\rm assemble}\, n_e\right).
\end{split}
\end{equation}
Therefore, $c_{\rm assemble} = 100$ and $c'_{\rm assemble} = 184$,  both constants are on the order of a few hundred.

The time step is chosen to satisfy both the advective CFL condition and the explicit viscous stability condition 
\[
\Delta t
\le
\min\left\{
\alpha_{\mathrm{CFL}}\frac{h}{\|v\|_{L^\infty}},
\frac{\alpha_{\mathrm{RK4}}}{\nu (2\pi k_{\max})^2}
\right\},
\]
here \(\alpha_{\mathrm{CFL}}\) is an \(O(1)\) CFL constant, $\|v\|_{L^\infty}$ is the maximum velocity magnitude, $\alpha_{\mathrm{RK4}}\approx 2.785$ is the stability limit of classical RK4~\cite{hairer1993solving}, and \(k_{\max}\) is the
largest retained Fourier index.
In our setting, the viscosity is small, so the advective CFL condition dominates. For fixed final time $T=1$ and \(\|v\|_{L^\infty}=O(1)\), this gives  $\Delta t = O(h)$ and $n_t = O(n)$. Thus,
\begin{equation}
C_{\rm solver}
\approx  c_{\rm assemble}\, n_e^{3/2}\log_2 n_e + c'_{\rm assemble}\, n_e^{3/2}.
\end{equation}

\subsection{Vehicle aerodynamics}
\label{supp-ssec:fluid}
The nonlinear steady RANS-\(k\)-\(\omega\) SST system consists of the three
momentum equations, the incompressibility constraint, and the two turbulence
transport equations for \(k\) and \(\omega\). In \texttt{OpenFOAM}, it is solved using a segregated
SIMPLE fixed-point iteration. At each SIMPLE iteration, \texttt{OpenFOAM} assembles
six separate linearized systems for the momentum, pressure-correction, \(k\), and
\(\omega\) equations, which are then solved sequentially in a
block-segregated manner. The momentum predictor is first solved using the
current pressure, turbulent viscosity, and face fluxes; the pressure-correction
equation is then solved to enforce discrete mass conservation and update the
velocity and flux fields. The corrected velocity and flux are subsequently used
to assemble and solve the \(k\)- and \(\omega\)-transport equations, after which
the eddy viscosity is updated for the next SIMPLE iteration.

In the \texttt{OpenFOAM} implementation, the implicit finite volume matrices are sparse and retain only the nearest-neighbor stencil induced by the mesh faces. 
For the convective
terms, only the lower-order part of the discretization is treated implicitly, while the
gradient-based linear-upwind correction is handled explicitly as a source contribution.
Similarly, for corrected diffusion schemes, the implicit matrix contains only the
uncorrected Gauss-Laplacian contribution (see \cref{supp-eq:corrected_diffusion}), while the non-orthogonal correction terms are
treated explicitly in the source term.
Let \(n_c\) denote the number of control volumes and \(n_f\) the number of
internal faces. Then each scalar finite volume matrix $A$ has approximately
\[
\operatorname{nnz}(A) = n_c + 2n_f,
\]
since each cell contributes one diagonal entry and each internal face
contributes two off-diagonal couplings. These matrices are reassembled at each SIMPLE iteration.

We estimate the cost of one SIMPLE iteration by separating it into an assembly cost and a
linear-solve cost. Consider first one representative scalar convection-diffusion
equation. The implicit matrix contains only the low-order convective contribution and the uncorrected
diffusive contribution. Since both are assembled by looping over the
internal faces, we estimate their costs as approximately \(4n_f\) and \(3n_f\)
flops, respectively.
The remaining assembly cost comes from explicit correction and source terms. A Gauss linear
gradient reconstruction for a cell-centered scalar $\phi_e$ is given by
 $$
 \nabla \phi_e \approx \frac{1}{|\Omega_e |}\sum_{f\subset\partial \Omega_e} \phi_f s_f,
 $$
 where \(|\Omega_e|\) is the cell volume, \(\phi_f\) is the face-interpolated value, and
\(s_f\) is the face area vector.
 We estimate this cost as approximately $12n_f + 3n_c$. 
 For the cell-limited Gauss gradient, the unlimited gradient is rescaled as
\begin{align}
 \widetilde{\nabla \phi_e} = \lambda_e  \nabla \phi_e  \qquad \phi_f \approx \phi_e + \widetilde{\nabla \phi_e}\cdot(x_f-x_e),   
\end{align}
 where \(0\le \lambda_e\le 1\) is chosen to prevent the reconstructed face values $\phi_f$ from creating new
extrema. We estimate the additional cost of limiter construction and gradient rescaling as $18n_f + 5n_c$.
Using these gradients, the explicit linear-upwind correction at face \(f\) is 
\begin{align}
\delta\phi_f^{\mathrm{up}}
=
\widetilde{\nabla \phi}_{e_{\mathrm{up}}}\cdot(x_f-x_{e_{\mathrm{up}}}),
\end{align}
and contributes \(F_f\,\delta\phi_f^{\mathrm{up}}\) to the source term, where \(F_f\) is
the face flux and \(e_{\mathrm{up}}\) is the upwind cell. We estimate this cost as $11n_f+n_c$.
For corrected diffusion, the diffusive face flux is decomposed into an implicit
orthogonal part and an explicit non-orthogonal correction:
\begin{equation}
\label{supp-eq:corrected_diffusion}
\begin{aligned}
    \nabla\cdot (\nu \nabla \phi_e) 
    &\approx \frac{1}{|\Omega_e|}\sum_{f\subset \partial \Omega_e} \nu_f (\nabla \phi)_f  \cdot s_f \\
    &\approx \frac{1}{|\Omega_e|}\sum_{f\subset \partial \Omega_e} \Bigl( 
    \nu_f \alpha_f (\phi_{e'} - \phi_{e})
    +
    \nu_f (s_f - \alpha_f d_{e,e^{'}}) (\nabla\phi)_{f}
    \Bigr),
\end{aligned} 
\end{equation}
where \(d_{e,e^{'}}=x_{e'}-x_e\) is the vector connecting the neighboring cell centers,
and \(\alpha_f\) is the scalar coefficient of the orthogonal two-point diffusive flux.
The second term is treated explicitly, and we estimate its cost as $16n_f+n_c$.
Therefore, for one representative scalar equation, the assembly cost is approximately
\begin{align}
 7n_f + (12n_f + 3n_c) + (18n_f + 5n_c) + (11n_f+n_c) +(16n_f+n_c) = 64n_f + 10n_c.
\end{align}
Counting the three velocity components, the pressure-correction equation, and the two
turbulence equations as six scalar-equivalent systems gives the rough assembly estimate
\begin{align}
C_{\rm assemble} 
\approx
6(64n_f+10n_c) = 384n_f+60n_c.
\end{align}
For the solve cost, each scalar finite volume matrix has approximately  
$n_c+2n_f$ nonzero entries. One symmetric Gauss-Seidel sweep therefore costs approximately \(4(n_c+2n_f)\) flops. If each scalar system requires \(m\)
inner iterations to reach the prescribed residual tolerance, and if all six systems are
modeled with the same per-iteration cost, then
\[
C_{\rm solve} \approx 6m\,4(n_c+2n_f) = 24m\,n_c + 48 m\, n_f.
\]
In practice, the pressure-correction equation is solved with the geometric agglomerated algebraic multigrid solver \texttt{GAMG}, while the velocity and turbulence equations are solved using \texttt{smoothSolver} with symmetric Gauss-Seidel smoothing. Their actual constants therefore differ, and the simplified model above should be understood only as a leading-order cost estimate.

Finally, if \(n_{t}\) SIMPLE iterations are required, and we use the approximation \(n_f\approx 3n_c\) for shape-regular finite volume meshes,  then the total cost is
approximately
\begin{align*}
C_{\rm solver}
\approx
n_{t}
(C_{\rm assemble} + C_{\rm solve}) 
=
n_{t}
(1212 n_c + 168 m \, n_c) 
= n_{t}
(c_{\rm assemble} n_c + c_{\rm solve} m \, n_c),
\end{align*}
where \(c_{\mathrm{assemble}} = 1212\) and \(c_{\mathrm{solve}} = 168\), which are on the
order of a few thousand and a few hundred, respectively. 
For the inner linear solves, we follow the settings of the official \texttt{OpenFOAM}~12 \texttt{drivaerFastback} case. The pressure, velocity, and turbulence equations use an absolute tolerance of \(10^{-6}\) and a relative tolerance of \(0.1\). Our solver logs show that each linear system typically requires two to three inner iterations. To avoid overstating the classical solver cost, we use \(m=2\) as a conservative effective value in the floating-point cost estimate.

\section{Additional numerical results}
\label{supp-sec:numerics}
This section provides additional numerical results that complement the
comparisons in the main paper. First, we test the sensitivity of the Darcy flow results to the neural operator architecture by repeating the experiment with a standard Fourier neural operator. Second, we examine the sensitivity of the Navier--Stokes results
to the choice of time integrator. Finally, we report the configuration and results of a representative \texttt{OpenFOAM} calculation for the canonical \texttt{drivaerFastback} geometry.

\subsection{Architecture sensitivity for Darcy flow}
\label{supp-sec:numerics-sensitivity}
\begin{figure}[htbp]
     \centering
     \includegraphics[width=0.9\textwidth]{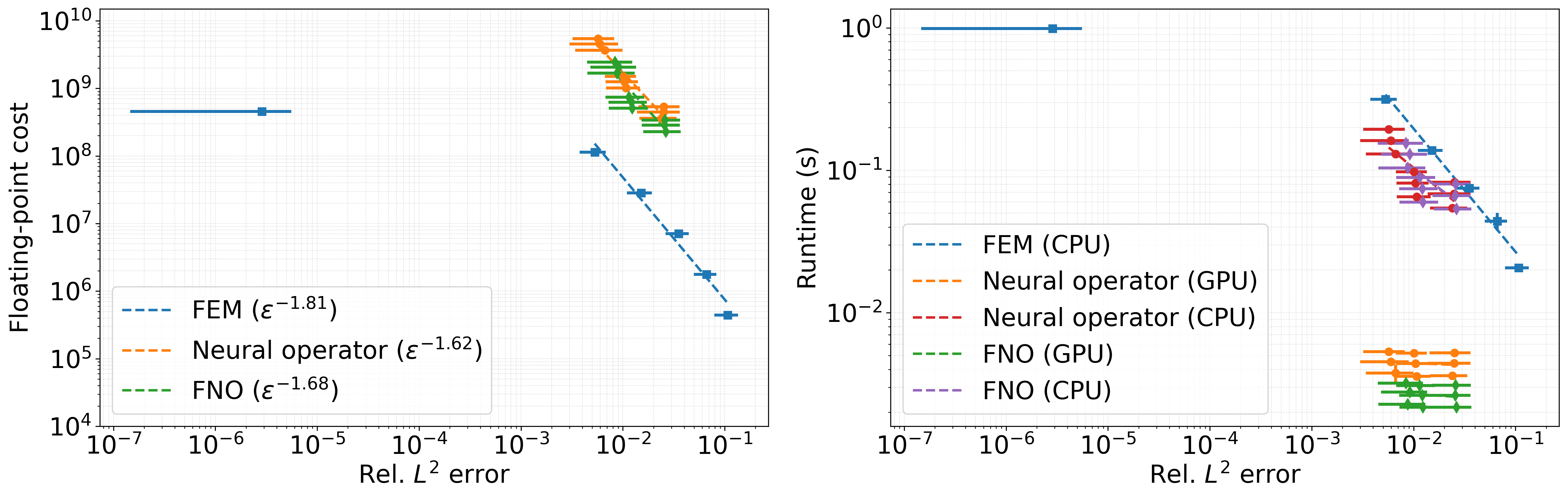}\\
     \caption{Darcy flow: cost--accuracy comparison among finite elements with geometric multigrid (FEM), the modified Fourier neural operator (Neural operator), and the original Fourier neural operator (FNO).
     Left: floating-point operation count versus relative $L^2$ error. 
     Right: wall-clock runtime in seconds versus relative $L^2$ error for CPU and GPU 
     implementations. Error bars denote one standard deviation over the test samples.}
     \label{supp-fig:darcy-cost-accuracy}
\end{figure}

To assess sensitivity to architecture, we repeat the Darcy flow experiment using a standard Fourier neural operator (FNO)~\cite{li2020fourier}, with the same data, training protocol, and evaluation procedure. As shown in
\cref{supp-fig:darcy-cost-accuracy}, the qualitative conclusion is unchanged: geometric multigrid requires less floating-point work than either neural operator, whereas both neural operators achieve shorter wall-clock runtimes over the accuracy range they attain, owing to the
efficient execution of dense tensor operations on modern hardware. These results suggest that the observed behavior is robust to this architectural modification. 
Because the underlying acceleration mechanisms are shared by
many neural operator architectures, we expect similar qualitative behavior more broadly; establishing this generality, however, requires comparisons across additional architectures and problems.

\subsection{Time integrator sensitivity for Navier-Stokes flow}
\label{supp-sec:numerics-etd}

To assess sensitivity to the choice of time integrator, we repeat the Navier-Stokes experiment using a first-order exponential
time differencing scheme (ETD1) and a fourth-order exponential time differencing Runge-Kutta scheme (ETDRK4)~\cite{cox2002exponential},
while retaining the same spatial discretization, data, and evaluation protocol. In our experiments, ETD1 is unstable at both the time step used for RK4 and ETDRK4 and at half that value. We therefore use
$\Delta t_{\rm ETD1}=\frac{\Delta t_{\rm RK4}}{4}=\frac{\Delta t_{\rm ETDRK4}}{4}$. ETD1 consequently takes four
times as many time steps but requires only one nonlinear right-hand-side
evaluation per step, whereas RK4 and ETDRK4 require four such evaluations.
The leading-order floating-point costs of the three methods are therefore
comparable under these settings.

\begin{figure}[htbp]
     \centering
     \includegraphics[width=0.9\textwidth]{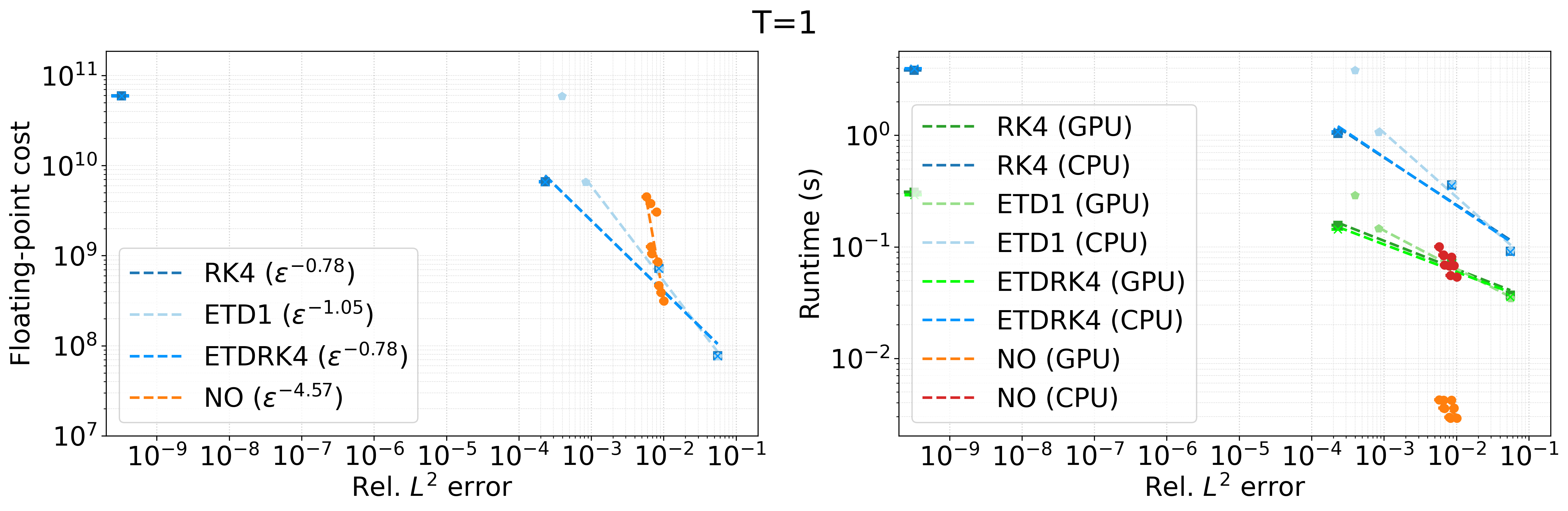}\\
     \includegraphics[width=0.9\textwidth]{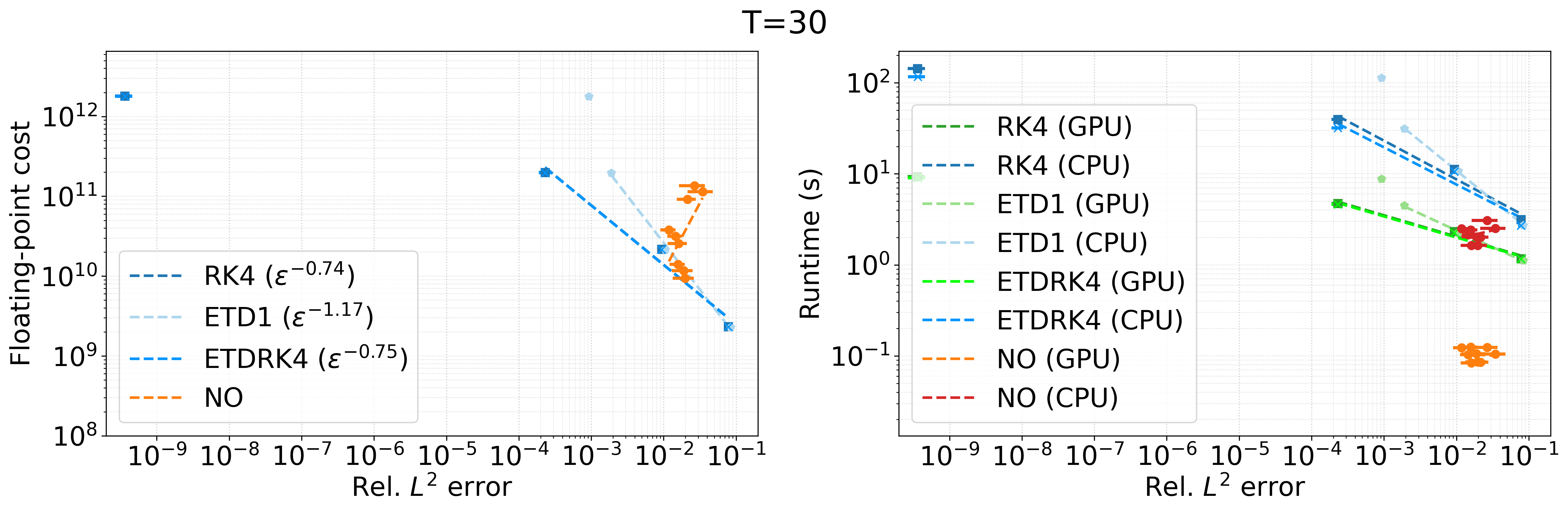}\\
     \caption{Time integrator sensitivity for incompressible
    Navier-Stokes flow: cost-accuracy comparisons among Fourier
    pseudospectral solvers using RK4, ETD1, and ETDRK4 and the neural
    operator (NO).  Top: teacher-forced one-step evaluation \((T=1)\), with each reference
    state \(\omega(t)\), \(t=0,\ldots,29\),used to
    predict \(\omega(t+1)\). Bottom: recurrent rollout from \(t=0\)
    through \(T=30\).
     Left: floating-point operation count versus relative $L^2$ error. 
     Right: wall-clock runtime in seconds versus relative $L^2$ error for CPU and GPU 
     implementations. Error bars denote one standard deviation over the test samples.}
     \label{supp-fig:NS-cost-accuracy}
\end{figure}

As shown in \cref{supp-fig:NS-cost-accuracy}, ETD1 yields larger trajectory prediction errors at comparable floating-point work, reflecting its first-order temporal accuracy and underscoring the importance of higher-order time integration. ETDRK4 and RK4 exhibit similar cost-accuracy behavior. In this convection dominated regime, exact treatment of the viscous linear term provides little additional benefit because the time step remains governed primarily by the explicitly treated nonlinear advection. Consequently, the qualitative comparison with the neural operator is unchanged by the choice among the time integrators considered here.

\subsection{Representative \texttt{OpenFOAM} vehicle configuration}
\label{supp-ssec:openfoam-fastback}
The vehicle aerodynamics results reported in the main paper are averaged over
the six test geometries. To illustrate the scale of an individual calculation
and facilitate reproducibility,
\cref{supp-tab:openfoam-fastback} reports the mesh sizes, iteration counts, CPU allocations, runtimes, and relative errors for the canonical
\texttt{drivaerFastback} geometry. These values are close to the corresponding test set averages reported in \cref{tab:openfoam_setup}. 

Although \texttt{OpenFOAM} runtimes and numerical errors vary across geometries, this variability can be monitored using mesh quality measures, residual histories, and the stabilization of aerodynamic force coefficients. These diagnostics make the sources of error in the classical solver more readily interpretable than those in the neural operator predictions. Moreover, because the classical solver is not trained on a finite dataset, it is not subject to out-of-distribution errors in the same statistical sense as the neural operator.

\begin{table}[htbp]
\centering
\begin{tabular}{l|cccc}
\Xhline{1.1pt}
Configuration
    & Large reference & Large & Medium & Small \\
\hline
Control volumes $n_c$                    & 22,551,532  & 22,551,532   & 2,982,335    & 440,405   
\\
SIMPLE iterations $n_t$&  7000  &  2000      & 1000 & 1000 
\\
CPU cores               & 256   & 256  & 32        & 4                  \\   
\hline
Runtime (s)              &   6424 &   1868    & 975     & 557                 \\
Relative $L^1$ error      & Reference &   11.46\%   & 16.31\%     & 23.62\%    \\
\Xhline{1.1pt}
\end{tabular}
\caption{
\texttt{OpenFOAM} configurations and results for prediction of surface pressure coefficient $C_p$ in the canonical \texttt{drivaerFastback} benchmark. 
Errors are measured relative to the reference solution computed on the large mesh using 7000 SIMPLE iterations.}
\label{supp-tab:openfoam-fastback}
\end{table}

\bibliographystyle{siamplain}
\bibliography{references}
\end{document}